\documentclass[aos,preprint]{imsart}
\RequirePackage[OT1]{fontenc}
\RequirePackage{amsthm,amsmath}
\RequirePackage[numbers]{natbib}
\RequirePackage[colorlinks,citecolor=blue,urlcolor=blue]{hyperref}
\usepackage{bm}
\startlocaldefs
\numberwithin{equation}{section}
\theoremstyle{plain}

\usepackage{amssymb}
\usepackage{amsthm,amssymb}
\usepackage{amsfonts}
\usepackage{numprint}
\npthousandsep{,}

\usepackage{graphicx}
\usepackage{verbatim,color,amssymb,epsfig}
\usepackage[usenames,dvipsnames]{xcolor}
\usepackage{fancyhdr}
\usepackage{alg}
\usepackage{natbib}

\usepackage{numprint}

\newcommand{\E}{\mathrm{E}}
\newcommand{\var}{\mathrm{Var}}
\newcommand{\Jn}{F_n}
\newcommand{\Ln}{S_{n,{\scriptscriptstyle SUB}}}
\newcommand{\Lni}{S_{n,{\scriptscriptstyle SUB}}^{{\scriptscriptstyle \mathrm{ID}}}}

\newcommand{\hatsigmaY}{\hat{\sigma}_{Y,n,{\scriptscriptstyle SUB}}^2}
\newcommand{\TnY}{T_{Y,n}}
\newcommand{\LnY}{S_{Y,n,{\scriptscriptstyle SUB}}}

\newcommand{\Lt}{\tilde{S}_{n}}

\newcommand{\lt}{\tilde{S}_{\ell}}
\newcommand{\Ll}{S_{\ell,{\scriptscriptstyle SUB}}}
\newcommand{\Lnt}{\tilde{S}_{n,{\scriptscriptstyle SUB}}}

\newcommand{\mn}{m_{n,{\scriptscriptstyle SUB}}}
\newcommand{\ml}{m_{\ell,{\scriptscriptstyle SUB}}}
\newcommand{\Cn}{C_{n,k_n}}
\newcommand{\Cnu}{\tilde{C}_{n,k_n}}

\newcommand{\Cniid}{C_{n,k_n}^{{\scriptscriptstyle \mathrm{ID}}}}

\newcommand{\Cniidb}{C_{n,k_nb}^{{\scriptscriptstyle \mathrm{ID},boot}}}
\newcommand{\Cl}{C_{\ell,k_\ell}}
\newcommand{\Clu}{\tilde{C}_{\ell,k_\ell}}
\newcommand{\hatsigma}{\hat{\sigma}_{n,{\scriptscriptstyle SUB}}^2}

\newcommand{\hatsigmaid}{(\hat{\sigma}_{n,{\scriptscriptstyle SUB}}^{{\scriptscriptstyle \mathrm{ID}}})^2}
\newcommand{\hatsigmal}{\hat{\sigma}_{\ell,{\scriptscriptstyle SUB}}^2}

\endlocaldefs
\begin{document}
\begin{frontmatter}
\title{Convolved subsampling estimation with
applications to block bootstrap}
\runtitle{Convolved subsampling}
\begin{aug}

\author{\fnms{Johannes} \snm{Tewes}\thanksref{t1}
\ead[label=e1]{johannes.tewes@ruhr-uni-bochum.de}},
\author{\fnms{Dimitris N.} \snm{Politis}\ead[label=e2]{dpolitis@ucsd.edu}}
\and\\
\author{\fnms{Daniel J.} \snm{Nordman}\thanksref{t2}
\ead[label=e3]{dnordman@iastate.edu}}
\thankstext{t1}{Research supported by by the Collaborative Research Center ``Statistical modeling of nonlinear dynamic processes'' (SFB 823, Teilprojekt C3) of the German Research Foundation (DFG).}
\thankstext{t2}{Research partially supported by NSF DMS-1406747.}

\runauthor{Tewes, Politis, Nordman}
\affiliation{Ruhr-Universit\"{a}t Bochum, University of California-San Diego, \\Iowa State University}
%
%NA 02/71
%
\address{Department of Mathematics\\
Ruhr-Universit\"{a}t Bochum\\
44780 Bochum, Germany\\
\printead{e1}\\
\phantom{E-mail:\ }}
\address{
Department of Mathematics\\
University of California, San Diego\\
La Jolla, CA 92093-0112, USA \\
\printead{e2}\\
\phantom{E-mail:\ }}
\address{Department of Statistics\\
Iowa State University\\
Ames, IA USA 50011\\
\printead{e3}\\
\phantom{E-mail:\ }}
\end{aug}
\begin{abstract}

The block bootstrap  approximates sampling distributions from dependent
data by resampling data blocks.  A fundamental
problem is establishing its consistency for   the distribution of a sample mean, as
a prototypical statistic. We use a structural relationship with
subsampling  to characterize the bootstrap
in a new and general manner. While subsampling and block bootstrap differ, the block
bootstrap distribution of a sample mean equals that of a $k$-fold self-convolution of a subsampling
distribution. Motivated by this, we provide simple necessary and sufficient conditions for a convolved
subsampling estimator to produce a normal limit that matches the target of bootstrap
estimation.  These conditions may be linked to consistency properties of an original subsampling distribution,
which are often obtainable under minimal assumptions. Through several examples, the results are shown to
validate the block bootstrap for  means  under significantly
weakened assumptions in many existing (and some new) dependence settings, which also addresses    a standing conjecture of  Politis, Romano and Wolf~(1999).
Beyond  sample means,  the convolved subsampling estimator  may not match the block
bootstrap, but instead provides a hybrid-resampling estimator of interest in its own right.  For general statistics with normal limits, results
also establish the consistency of convolved subsampling under minimal dependence conditions, including non-stationarity.
\end{abstract}

\end{frontmatter}

\section{Introduction}

Subsampling and  block bootstrap  are two   common  nonparametric tools for statistical inference under dependence;
 see  \citet{PoRoWo} and \citet{Lah2}, respectively, for monographs on these. Both aim to approximate   distributions of  statistics with correlated data, and both are data resampling methods that
 use blocks of neighboring observations to capture   dependence.
  The subsampling approach of  \citet{PoRo1} treats data blocks as  small scale renditions of the original data, which   provides
  replication of a statistic for estimating a sampling distribution.
  The block bootstrap  differs philosophically by using data blocks as building material
  to re-create the original data.  Essentially,  the block bootstrap selects  and glues blocks  to reproduce a full-scale set of data, as proposed by \citet{Kuen2} and \citet{LiSi}  for extending \citet{Efr}'s
  bootstrap to time series.   As noted in \citet{PoRoWo} (cf.~sec.~3.9), subsampling is often valid under weak assumptions about
  the dependent process,  basically requiring that a non-degenerate (possibly non-normal) limit exist for the sampling distribution being approximated.  In contrast, the block bootstrap applies to mean-like statistics with normal limits and typically requires comparatively much stronger assumptions for its validity.   Case-by-case treatments are commonly needed to validate the bootstrap across differing dependence conditions.  However, while perhaps not widely recognized,   subsampling can in fact be used to verify the block bootstrap in some cases, which is a theme of this work.

We investigate estimators defined by the $k$-fold self-convolution of a  subsampling distribution,
and establish a new and general theory for their consistency to normal limits.  There are two basic motivations for considering such convolved subsampling.   The first is that,
in the fundamental case of sample means, the block bootstrap estimator is  a $k$-fold self-convolution of a subsampling distribution  (centered and normalized), where the level $k$ of convolution corresponds to the number of resampled blocks. This observation
 was originally noted by \cite{PoRoWo}, who suggested   this aspect as a potential technique for showing the validity of the bootstrap.
 Specifically, they conjectured that convolved subsampling might provide a route for establishing the block bootstrap
 under minimal conditions for non-stationary, strongly mixing processes, in analogy to bootstrap results
 existing for   stationary, mixing series due to \citet{Radu1996, Radu2012}.  For the bootstrap under dependence, the
  findings of \citet{Radu1996, Radu2012} for the sample mean have stood out as an exception, verifying the method under the same weak assumptions as subsampling (i.e., conditions essentially needed for a limit law to exist).
By investigating the convolved subsampling approach here,
 we can answer the above conjecture   affirmatively.
Moreover, we show convolved subsampling leads to a simple and unified procedure for establishing the block bootstrap for sample means
 for further types of processes under much weaker conditions than previously considered,  such as linear time processes, long-memory sequences,  (non-stationary) almost periodic time series, and
 spatial fields.  Hence, convolved subsampling estimation allows for  bootstrap consistency under dependence to be generally extended under the same weak
 assumptions used by subsampling,
containing the conclusions of \citet{Radu1996, Radu2012} for stationary  series as a special case.

  While connections to the bootstrap are useful,  our study of  convolved subsampling estimation
  is intended  to be broad, applying also to
general statistics with normal limits and with arbitrary levels of convolution.  Consistency results often do not require particular assumptions about the underlying
dependent process, but are rather formulated in terms of
   mild convergence properties of the original subsampling distribution and its variance.   Furthermore, we show that a consistent subsampling variance is not only sufficient, but essentially necessary, for the consistency of convolved subsampling  (and the block bootstrap in some cases).  Due to its importance, we also provide tools for verifying the consistency of subsampling variance estimators.

For general statistics beyond the sample mean, the convolved subsampling distribution may differ from the block bootstrap,
which relates to a second motivation for our development.  That is, a general theory for convolved subsampling is of interest in its own right, as the approach can
 be  computationally less demanding than the block bootstrap while also potentially enhancing ordinary subsampling
for approximating sampling distributions with normal limits.  In fact, there has been recent interest
in establishing generalized types of  subsampling estimation for complicated statistics under various dependence structures,
where numerical studies suggest such
methods exhibit better finite sample performance than standard subsampling; for example, see  \citet{Len2} and \citet{ShTeWe2} for spectral estimates and U-statistics, respectively, with time series.
 While not formally recognized as such, however,
these proposed methods are exactly convolved subsampling estimators.  By exploiting this realization, our results can
facilitate future work and allow such previous findings with generalized subsampling to be demonstrated in an alternative, simpler manner with weaker assumptions.

 Section~\ref{sec:prelim} describes convolved subsampling estimation
and its connection to the block bootstrap.     General distributional results
for convolved subsampling are given in Section~\ref{sec:main}, while Section~\ref{sec:app} presents some applications
  with differing dependence structures.     Section~\ref{sec:mix} provides a broad result for
  convolved subsampling estimation of statistics from mixing time series.
  Under weak conditions, Sections~\ref{sec:mix2}-\ref{sec:spatial} apply convolved subsampling for demonstrating the block bootstrap for sample means
  with non-stationary time series (Section~\ref{sec:mix2} and the conjecture of \citet{PoRoWo}), linear time processes
  (Section~\ref{sec:linear}), long-range dependence (Section~\ref{sec:lrd}), and spatial data (Section~\ref{sec:spatial}).
   Section~\ref{sec:convmore} then briefly describes some relationships to other recent work with generalized subsampling, and
Section~\ref{sec:indep}    provides a short treatment of   independent data.
   Concluding remarks are given in Section~\ref{sec:concl}, and the proofs of main results appear in a supplement  \citep{TPN17}.
\section{Description of  convolved subsampling estimators}
\label{sec:prelim}
\subsection{Problem background and original subsampling estimation}
\label{sec:sub}
Consider data $X_1,\ldots,X_n$ from a real-valued process governed by a probability structure $P$.  For concreteness,
we may envision such observations arising from a   time series process $\{X_t\}$,
though spatial and other data schemes may be treated as well.
Based on  $X_1,\ldots,X_n$, consider the problem of approximating the distribution of
\[T_n\equiv \tau_n (t_n(X_1,\ldots,X_n) - t(P)),\] involving an estimator $t_n\equiv t_n(X_1,\ldots,X_n)$ of a parameter $t(P)$
and a sequence of positive scaling factors $\tau_n$ yielding a distributional limit for $T_n$.  For example, if
$t_n(X_1,\ldots,X_n) \equiv \bar{X}_n=\sum_{i=1}^n X_i/n$ is the sample mean, then $t(P)$ may correspond to a common process mean $\mu$
and $T_n$ may be defined with usual scaling $\tau_n=\sqrt{n}$   under weak time dependence.
Denote the sampling distribution function of $T_n$ as $\Jn(x) = P(T_n \leq x)$, $x\in \mathbb{R}$.

We next define the subsampling estimator of $\Jn$; see \citep{PoRo1}.  For  a positive integer $b\equiv b_n <n$, let $\{ (X_i,\ldots, X_{i+b-1}): i=1,\ldots,N_n\}$ denote the set of  $N_n \equiv n-b+1$  overlapping  data blocks, or subsamples, of length $b$.  To keep blocks  relatively small,
  the block size is often assumed to satisfy $b^{-1}+b/n +\tau_b/\tau_n\rightarrow 0 $ as $n\to \infty$.
For each subsample,
we compute the statistic as $t_{n,b,i} =t_b(X_{i},\ldots,X_{i-b+1})$ and  define a ``scale $b$" version of $T_n\equiv \tau_n (t_n(X_1,\ldots,X_n) - t(P))$ as
        $\tau_b [t_{n,b,i} -  t_n ]$ for $i=1,\ldots,N_n$.  Letting $I(\cdot)$ denote the indicator function,
  the subsampling estimator of $\Jn$ is given by
 \begin{equation}
 \label{eqn:sub}
 \Ln(x) = \frac{1}{N_n}\sum_{i=1}^{N_n} I\big( \tau_b [t_{n,b,i} -  t_n ] \leq x\big), \quad x\in \mathbb{R},
 \end{equation}
or the empirical distribution of subsample analogs $\{\tau_b [t_{n,b,i} -  t_n ]\}_{i=1}^{N_n}$ (cf.~\cite{PoRoWo}).

Suppose that  $\Ln$ is consistent for the distribution of $T_n$, where the latter has an asymptotically normal $N(0,\sigma^2)$ limit for some $\sigma>0$, that is, as $n\to \infty$,
\begin{equation}
\label{eqn:1}
    T_n  \stackrel{d}{\rightarrow}  N(0,\sigma^2),
    \end{equation}
    \begin{equation}
    \label{eqn:2}
       \sup_{x \in \mathbb{R}}|\Ln(x) - \Phi(x/\sigma) | \stackrel{p}{\rightarrow}  0,
\end{equation}
 where $\Phi(\cdot)$ is the standard normal distribution function.  We wish to consider
 estimators of the distribution $F_n$ of $T_n$ formed by
  self-convolutions of the
 subsampling estimator $\Ln$.    This provides a general class of block resampling estimators in its own right, but
 also has explicit connections to block bootstrap estimators in the important case that the statistic of interest $t_n(X_1,\ldots,X_n)=\bar{X}_n$ is a sample mean, as described next.
\subsection{Convolved subsampling and connections to block bootstrap}
\label{sec:conv}
 Let $k_n \in \mathbb{N}$ be a sequence of positive integers and define a triangular array $\{ Y_{n,1}^*,\ldots,Y_{n,k_n}^*  \}_{n \geq 1}$, where, for each $n$, $\{Y_{n,i}^* \}_{i=1}^{k_n}$ are iid variables following the subsampling distribution $\Ln$, as determined by (\ref{eqn:sub}) from data $X_1,\ldots,X_n$.  For $n\geq1$, define a centered and scaled sum
 \begin{equation}
 \label{eqn:zn}
   Z_n^* \equiv \frac{1}{\sqrt{k_n}}\sum_{j=1}^{k_n} (Y_{n,i}^* -\mn)
 \end{equation}
  where  $\mn \equiv \int x d \Ln(x) = N_n^{-1} \sum_{i=1}^{N_n}   \tau_b [t_{n,b,i} -  t_n ]   $ is the mean of the subsampling distribution $\Ln$, and let \[\Cn(x) \equiv P_*( Z_n^* \leq x), \quad x\in \mathbb{R},\] denote the induced resampling distribution $P_*$ of $Z_n^*$.  Then, $\Cn$ represents the $k_n$-fold self-convolution
  of the subsampling distribution $\Ln$, with appropriate centering/scaling adjustments. That is,
  \[
       \Cn(x) =  \underbrace{\Ln *  \Ln * \cdots * \Ln}_{\mbox{$k_n$ times}}(x \sqrt{k_n}   + k_n m_n ), \quad x\in \mathbb{R}.
  \]
  We consider $\Cn$ as an estimator of the distribution $\Jn$ of $T_n$
  and formulate general conditions under which this convolved subsampling distribution  is also consistent.

 As suggested earlier, such results have direct implications for block bootstrap estimation as well, because
the convolved subsampling estimator $\Cn$ exactly matches a block bootstrap estimator
in the basic  sample mean case
   $t_n(X_1,\ldots,X_n)=\bar{X}_n$.
 To illustrate, consider approximating the distribution of $T_n = \sqrt{n}(\bar{X}_n - \mu)$  where $t(P)\equiv \mu=\E \bar{X}_n$ and $\tau_n=\sqrt{n}$.
 In this setting, the block bootstrap uses an analog
 \begin{equation}\label{eqn:bb} T_n^* = \sqrt {n_1 } (\bar{X}^*_{n_1 } - \E_* \bar{X}^*_{n_1 })\end{equation} based on the average  $\bar{X}^*_{n_1}\equiv n_1^{-1}\sum_{i=1}^{n_1}X_i^*$
 from a block bootstrap sample $X_1^*,\ldots,X_{n_1}^*$ of size $n_1 \equiv k_n b$,   which is  defined by drawing $k_n$ blocks of length $b$, independently and with replacement, from the subsample collection $\{(X_i,\ldots,X_{i+b-1}):i=1,\ldots,N_n\}$
  and pasting these together (where above $\E_* \bar{X}^*_{n_1} = N_n^{-1}\sum_{i=1}^{N_n} b^{-1}\sum_{j=i}^{i+b-1} X_j$ denotes
   the bootstrap expectation of $\bar{X}^*_{n_1}$); see ch.~2, \citet{Lah2}. Most typically, the number of blocks resampled
is taken as $k_n = \lfloor n/b \rfloor \rightarrow \infty$ so that the bootstrap sample re-creates the approximate
length $ \lfloor n/b\rfloor b \approx n$ of the original sample.  The bootstrap distribution of $T_n^*$ here  is then equivalent to the convolved subsampling distribution $\Cn$. This is because
   $T_n^*$
 has the same
resampling distribution as $Z_n^*$
in (\ref{eqn:zn}) as a sum of  $k_n$ iid block averages $(Y_{n,i}^*-\mn)/\sqrt{k_n}$, with each $Y_{n,i}^*$ drawn
   from $\Ln$ in (\ref{eqn:sub}) where $t_n = \bar{X}_n$ and  $\tau_b[ t_{n,b,i}-t_n] = \sqrt{b}[b^{-1}\sum_{j=i}^{i+b-1} X_j-\bar{X}_n]$, $ 1 \leq i \leq  N_n$,  for the sample mean case.
     Consequently, if convolved subsampling estimators
$\Cn$ are shown to be valid under weak conditions, such results entail that block
bootstrap estimation is as well. In the following, we make comprehensive use of the fact that
$\Cn$ is always and exactly a block bootstrap estimator whenever the underlying statistic
$t_n(X_1, \ldots, X_n) = \bar{X}_n$
  is a sample mean; this holds true across all the varied dependent
data structures considered here, including settings where the block bootstrap formulation
(\ref{eqn:bb}) itself requires some modification (cf.~long-range dependence in Section~\ref{sec:linear}).

\section{Fundamental results for convolved subsampling}
\label{sec:main} From (\ref{eqn:sub}) and  the  subsampling  mean $ \mn\equiv \int  x  d \Ln(x) = N_n^{-1}\sum_{j=1}^{N_n} \tau_b[t_{n,b,i}-t_n]$,
  we  have the variance of the original subsampling distribution
 $\Ln$ as   \[\hatsigma\equiv \int (x-\mn)^2 d \Ln(x)= \frac{1}{N_n}\sum_{j=1}^{N_n} \left(\tau_b[t_{n,b,i}-t_n]-\mn\right)^2,\]  which estimates the asymptotic variance $\sigma^2$
 of $T_n$ as in (\ref{eqn:1}) (cf.~\citep{PoRoWo}).  Note that $\hatsigma$ is also the variance of the convolved subsampling distribution
$\Cn$ (i.e., the variance of the iid sum  from (\ref{eqn:zn})). Correspondingly, $\hatsigma$ is
then a block bootstrap variance estimator when applied to sample means.

Sections~\ref{sec:main1}-\ref{sec:main3} provide basic distributional results for convolved subsampling estimators,  describing when and how these have normal limits.  These findings do not involve
  particular assumptions about the  process $\{X_t\}$, but are instead  expressed
   through properties of the original subsampling distribution $\Ln$ and, specifically,  convergence of the subsampling variance $\hatsigma$.  Such subsampling properties can
  often be verified under weak assumptions about a process,  allowing the limit behavior of convolved estimators $\Cn$, and the block bootstrap, to be established under minimal conditions.  Results in Section~\ref{sec:main1} address the important case
  where the original subsampling distribution $\Ln$ has a normal limit (\ref{eqn:2}), as is often natural when
  the  statistic $T_n \stackrel{d}{\rightarrow} N(0,\sigma^2)$ is asymptotically normal. These findings are expected to be the most practical for establishing convolved subsampling  $\Cn$ estimation with normal  targets   (\ref{eqn:1}).   Dropping the condition that $\Ln$ converges to a  normal law but
  assuming  convolved estimators $\Cn$ are based on increasing convolution $k_n\to \infty$ of $\Ln$,
   Section~\ref{sec:main2} characterizes
  the convergence of $\Cn$ to normal limits through the subsampling variance $\hatsigma$.  In many problems involving the block bootstrap for sample means (cf.~Section~\ref{sec:app}), where $T_n$ has a normal limit (\ref{eqn:1}), these results provide both necessary and sufficient conditions for
  the validity of the block bootstrap as well as convolved subsampling generally.
       Finally, because  convergence  $\hatsigma \stackrel{p}{\rightarrow}\sigma^2$ of the  subsampling variance emerges as central to the behavior of
  convolved   estimators $\Cn$, Section~\ref{sec:main3} develops basic  results for establishing this feature.

\subsection{Convolution of subsampling distributions with normal limits}

\label{sec:main1}

Theorem~\ref{theorem1} provides a sufficient condition for the general validity of the convolved estimator
$\Cn$ via fundamental subsampling quantities, $\Ln$ and $\hatsigma$.
%  When the subsampling
%distribution $\Ln$ converges to a normal limit (\ref{eqn:2}) and the subsampling variance $\hatsigma$
%is also consistent, then {\it any} amount $k_n$
%of convolution of $\Ln$ is guaranteed to produce an estimator $\Cn$ with a normal limit, which will be then %consistent for the distribution $F_n$ of
%a quantity $T_n$ with a normal limit (\ref{eqn:1}).

\newtheorem{theoremnew1}{Theorem}

\begin{theoremnew1}
\label{theorem1} Suppose (\ref{eqn:2}) holds (i.e., $\sup_{x \in \mathbb{R}} |\Ln(x) - \Phi(x/\sigma)| \stackrel{p}{\rightarrow} 0$) and $\hatsigma \stackrel{p}{\rightarrow} \sigma^2>0$
as $n \to \infty$. Then,
\[
  \sup_{x \in \mathbb{R}} |\Cn(x) - \Phi(x/\sigma)| \stackrel{p}{\rightarrow} 0 \quad \mbox{as $n\rightarrow \infty$}
\]
for  any positive integer  sequence $k_n$.\\
Furthermore, when (\ref{eqn:1}) holds additionally (i.e., $T_n \stackrel{d}{\rightarrow} N(0,\sigma^2)$), then $\Cn$ is consistent for the distribution $\Jn$ of $T_n$,
\[
  \sup_{x \in \mathbb{R}} |\Cn(x) - \Jn(x)| \stackrel{p}{\rightarrow} 0 \quad \mbox{as $n\rightarrow \infty$}.
\]

\end{theoremnew1}

  To re-iterate, the integer  sequence $k_n$, $n \geq 1$, need not even be convergent in Theorem~\ref{theorem1}.
 The consistency of the subsampling variance estimator $\hatsigma$  automatically guarantees that, for {\it any} amount $k_n$
of convolution of $\Ln$, the convolved subsampling estimator $\Cn$ will have a normal limit if the subsampling distribution  $\Ln$ does. In other words, if
(\ref{eqn:1})-(\ref{eqn:2}) hold so that $\Ln$ is consistent, then $\Cn$ will be as well provided $\hatsigma \stackrel{p}{\rightarrow} \sigma^2$.
When the statistic $t_n(X_1,\ldots,X_n) = \bar{X}_n$
  is a sample mean, then $\Cn$ again denotes a
block bootstrap estimator based on $k_n$ resampled blocks, which is thereby consistent under
Theorem~\ref{theorem1} for any sequence $k_n$, including the common choice $k_n = \lfloor n/b \rfloor \rightarrow \infty$.

Proposition~\ref{prop1} next characterizes the convolved subsampling estimator $\Cn$ under {\it bounded} levels $k_n$ of  convolution.  In this case,  a normal limit for the subsampling estimator $\Ln$ entails the same for  the convolved  estimator $\Cn$, provided
the mean $\mn \equiv \int x d \Ln(x)$ of the subsampling distribution converges to zero.  But, if the subsampling mean $\mn$
converges in this fashion, a normal limit for $\Cn$ with bounded $\{k_n\}$ is equivalent to a normal limit for the original subsampling distribution $
\Ln$.

\newtheorem{proposition1}{Proposition}

\begin{proposition1}
\label{prop1} Suppose $\sup_n k_n < \infty$.\\
(i)   If  (\ref{eqn:2}) holds (i.e., $\sup_{x \in \mathbb{R}} |  \Ln(x) - \Phi(x/\sigma)| \stackrel{p}{\rightarrow} 0$),
then
\[
\sup_{x \in \mathbb{R}} |  \Cn(x) - \Phi(x/\sigma)| \stackrel{p}{\rightarrow} 0 \quad \mbox{as $n\to\infty$}
\]
if and only if  $\mn \equiv \int x d \Ln(x) \stackrel{p}{\rightarrow} 0 $.\\
(ii) If $\mn  \stackrel{p}{\rightarrow} 0 $  as $n\to \infty$, then  (\ref{eqn:2}) holds if and only if
\[
\sup_{x \in \mathbb{R}} |  \Cn(x) - \Phi(x/\sigma)| \stackrel{p}{\rightarrow} 0 \quad \mbox{as $n\to\infty$}.
\]

\end{proposition1}

When the original subsampling estimator $\Ln$ is consistent for a distribution with a
normal limit (i.e., (\ref{eqn:1})-(\ref{eqn:2})), both Theorem~\ref{theorem1} and Proposition~\ref{prop1} show that the convolved subsampling
estimator $\Cn$  is consistent under an additional subsampling moment condition.
With bounded levels $k_n$ of convolution, the additional condition under Proposition~\ref{prop1}  is
that the subsampling mean converge $\mn  \stackrel{p}{\rightarrow} 0 $. But, for general and potentially unbounded
$k_n$, the additional condition from Theorem~\ref{theorem1} for consistency of $\Cn$ is a convergent
subsampling variance $\hatsigma\stackrel{p}{\rightarrow}\sigma^2$. With  diverging amounts $k_n\to \infty$ of convolution,
which is often encountered in practice and in connection to the block bootstrap,
it turns out that convergence $\hatsigma\stackrel{p}{\rightarrow}\sigma^2$ is also {\it necessary} for consistency of the convolved
estimator $\Cn$, as treated in the next section.

\subsection{Unbounded convolution of subsampling distributions}
\label{sec:main2}
We next consider the behavior of convolved subsampling estimators for unbounded  convolution $k_n \to \infty$ as $n\to \infty$, which arises, for example, with the block bootstrap   $\Cn$
for sample means with $k_n = \lfloor n/b \rfloor $   resampled blocks.  Results here do not explicitly require convergence of the original subsampling estimator $\Ln$
 to a normal limit (\ref{eqn:2}).  While a reasonable condition in problems where the target quantity $T_n \stackrel{d}{\to} N(0,\sigma)^2$ is asymptotically normal, limits for $\Ln$ are not directly necessary
for convolved estimators $\Cn$ to yield valid normal limits from increasing convolution $k_n$
of $\Ln$.   However, convergence of
the subsampling variance $\hatsigma$ is   crucial, as shown in Theorem~\ref{theorem2}.

\begin{theoremnew1}
 \label{theorem2} Suppose $  k_n \to \infty $ and $ \int_{|x| \geq  \sqrt{k_n} \epsilon} x^2 d \Ln(x)  \stackrel{p}{\rightarrow}0 $ for each $\epsilon>0$ as $n \to \infty$. \\
(i)  Then,
  \begin{align*}
\sup_{x \in \mathbb{R}} |\Cn(x) - \Phi(x/\sigma)| \stackrel{p}{\rightarrow} 0
\end{align*}
 if and only if $\hatsigma\stackrel{p}{\rightarrow} \sigma^2>0$ as $n\rightarrow \infty$.\\
 (ii) When $\hatsigma\stackrel{p}{\rightarrow} \sigma^2>0$ as $n\rightarrow \infty$, then
  \begin{align*}
\sup_{x \in \mathbb{R}} |\Cn(x) - \Jn(x)| \stackrel{p}{\rightarrow} 0
\end{align*}
if and only if (\ref{eqn:1}) holds (i.e., $\sup_{x \in \mathbb{R}} |\Jn(x) - \Phi(x/\sigma)| \rightarrow  0$ for the distribution $\Jn$ of $T_n$).
\end{theoremnew1}

For an {\it unbounded} sequence $k_n\to\infty$ of convolution  (e.g., block bootstrap with $k_n = \lfloor n/b \rfloor $
concatenated blocks), Theorem~\ref{theorem2} imposes no direct assumption on the convergence of the
original subsampling distribution, but rather that $\Ln$ fulfills a mild
truncated second moment property. From this, the consistency of the convolved subsampling
estimator $\Cn$ to a normal limit is completely be determined by the behavior of the
subsampling variance $\hatsigma$ under Theorem~\ref{theorem2}. Furthermore, when $\hatsigma$ converges, the
convolved estimator $\Cn$ is again valid for estimating the distribution $\Jn$ of the target
quantity $T_n$ with a normal limit (Theorem~\ref{theorem2}(ii)).

The following corollary  of Theorem~\ref{theorem2} shows that a convolved estimator
  $\Cn$ will quite generally have a normal limit,
     provided that the subsampling variance converges $\hatsigma\stackrel{p}{\rightarrow} \sigma^2>0$  and that some other
     basic feature exists for the subsampling distribution $\Ln$ or for composite statistics $\{\tau_b [t_{n,b,i} - t_n] \equiv \tau_b [t_b(X_i,\ldots,X_{i-b+1}) - t_n(X_1,\ldots,X_n) ] \}_{i=1}^{N_n \equiv n-b+1} $  defining $\Ln$ in (\ref{eqn:sub}).
 Essentially, Corollary~\ref{cor1} entails that the truncated second moment assumption in Theorem~\ref{theorem2} is mild in conjunction with  $\hatsigma\stackrel{p}{\rightarrow} \sigma^2$.

\newtheorem{corollary1}{Corollary}
\begin{corollary1}
\label{cor1} Suppose one of the following conditions (C.1)-(C.4) holds:\\
(C.1)   for some distribution $J_0$ with variance $\sigma^2>0$, $\Ln(x) \stackrel{p}{\rightarrow} J_0(x)$ as $n\to \infty$   for any continuity point $x \in \mathbb{R}$ of $J_0$; \\
(C.2)  for some $\epsilon_0>0$, $N_n^{-1} \sum_{i=1}^{N_n} \big[\tau_b (t_{n,b,i} -t_n) \big]^{2+\epsilon_0}=O_p(1)$.\\
(C.3) the subsample-based sequence $\{T_{b,i}^2 \equiv \tau_b^2[t_{n,b,i} -t(P)]^2 : i =
1,\ldots, N_n\}_{n\geq 1}$ is uniformly integrable and $T_n \equiv  \tau_n(t_n- t(P)) = O_p(\tau_n/\tau_b)$\\
(C.4) $\{X_t\}$ is stationary, $\{T_n^2 : n \geq 1\}$ is uniformly integrable, and $\tau_b/\tau_n = O(1)$.\\[.1cm]
 Then,  as $n\to \infty$,
\[
\sup_{x \in \mathbb{R}} |  \Cn(x) - \Phi(x/\sigma)| \stackrel{p}{\rightarrow} 0
\]
for any sequence $k_n  $ with $\lim_{n\to \infty} k_n=\infty$ if and only if  $\hatsigma\stackrel{p}{\rightarrow} \sigma^2>0$
\end{corollary1}
\noindent {\sc Remark 1}: For reference, note $\tau_b/\tau_n \to 0$ often holds with subsample scaling as $n\to \infty$  so that, for example, the conditions
 $\tau_b/\tau_n=O(1)$ and $T_n   = O_p(\tau_n/\tau_b)$ are mild.\\

Hence, if $k_n\to\infty$ and  $\hatsigma \stackrel{p}{\rightarrow}\sigma^2$, then the convolved estimator $\Cn$ will converge to a normal limit if the subsampling distribution $\Ln$ is convergent (C.1) or  has an appropriate stochastically bounded  moment  (C.2), or if
the subsampling statistics related to computing $\Ln$ have uniformly integrable second moments (C.3)-(C.4).  Condition (C.4)
is a special case of (C.3) under stationarity, and corresponds to an underlying assumption of \citet{Radu1996, Radu2012}
for examining the block bootstrap estimator $\Cn$ of a sample mean with stationary, mixing processes; see also Remark~2 to follow.
When restricted to Condition (C.1), the ``$\Leftarrow$" part of Corollary~\ref{cor1} corresponds to an initial
convolved subsampling result due to \citet{PoRoWo} (Proposition~4.4.1) for unbounded convolution $k_n\to \infty$,  which was developed for
 establishing the block bootstrap estimator $\Cn$ for the sample mean of non-stationary data, as re-considered here in Section~\ref{sec:mix2}.
 Note that, for   inference with  $T_n$ having a normal $N(0,\sigma^2)$ limit (\ref{eqn:1}), Condition C.1 in Corollary~\ref{cor1} is perhaps most natural
 and approached by verifying convergence $\Ln$ to a normal (\ref{eqn:2}).     In which case, the implication of Corollary~\ref{cor1} (involving $k_n\to \infty$) for guaranteeing that convolved subsampling and  block bootstrap estimators replicate normal limits when $\hatsigma \stackrel{p}{\to}\sigma^2$ also becomes a special case of Theorem~\ref{theorem1} (involving {\it any} $k_n$).\\

\noindent {\sc Remark 2}: For block bootstrap estimation of the sample mean $T_n =
\sqrt{n}( \bar{X}_n-\E X_1)$
  with
strongly mixing, stationary processes, \citet{Radu1996, Radu2012}  provides necessary and sufficient conditions for convergence of $\Cn$ (with $k_n = \lfloor n/b \rfloor \to \infty$) to a normal limit, assuming
$\{T_n^2 : n \geq 1\}$ is uniformly integrable. Under such assumptions, the main result there is that normal
limits for both $\Cn$ and $T_n$ are equivalent. In comparison, the necessary and sufficient
conditions for normality of the block bootstrap estimator $\Cn$ for a mean in Theorem~\ref{theorem2} are perhaps
more basic in that the conclusions of  \citep{Radu1996, Radu2012}, under the additional
assumptions made there, follow from Theorem~\ref{theorem2} (cf.~Corollary~\ref{cor1}). In this sense, Theorem~\ref{theorem2} broadly re-frames  the findings in \citep{Radu1996, Radu2012}, by not involving particular process
assumptions (i.e., stationarity or mixing) and applying to convolved subsampling estimators
$\Cn$ with general statistics and arbitrarily increasing convolution levels $k_n \to \infty$. Further
connections to, and extensions of, the   results of \citet{Radu1996, Radu2012} are
made in Section~\ref{sec:mix} for strongly mixing processes.

 \subsection{Consistency of subsampling variance estimators }
\label{sec:main3}
Theorems~\ref{theorem1}-\ref{theorem2} demonstrate that the subsampling variance  $\hatsigma$ plays a key
role in the convergence of the convolved subsampling estimator $\Cn$ generally, and of the block
bootstrap for the sample mean in particular. However, convergence of the subsampling
distribution $\Ln$ itself is often much easier to directly establish under weak assumptions about
the process $\{X_t\}$; see \citet{PoRoWo} and Section~\ref{sec:app} to follow. This raises a further
question considered next: if one knows that subsampling estimator $\Ln$  is consistent (\ref{eqn:2})
for a normal limit, then when will the subsampling variance  $\hatsigma$ be convergent
as well, thereby guaranteeing (from Theorem~\ref{theorem1}) that the convolved estimator
$\Cn$ is also consistent? As shown in Theorem~\ref{theorem3}, a general characterization is possible
as well as simple sufficient conditions based on moment properties of  subsample statistics (e.g., $T_b^2$).

For $n \geq 1$, recall $T_n \equiv \tau_n(t_n(X_1,\ldots,X_n) - t(P))$ and additionally define $T_{n,i} \equiv
\tau_n(t_n(X_i,\ldots, X_{i+n-1})- t(P))$ for $i \geq 1$ from the statistic applied to $(X_i,\ldots,X_{i+n-1})$.
Based on $N_n \equiv n-b+1$  subsample observations of length  $1 \leq b \equiv b_n   < n$, define a
distribution function
\begin{equation}
\label{eqn:Dn}
D_{n,b}(x) \equiv \frac{1}{N_n} \sum_{i=1}^{N_n} P(T_{b,i}\leq x),\quad x\in\mathbb{R},
\end{equation}
as an average of subsample-based probabilities.

\begin{theoremnew1}  Suppose (\ref{eqn:2}) and $T_n = o_p(\tau_n/\tau_b)$ as $n \to \infty$.
\label{theorem3} \\
(i) Then, $\hatsigma\stackrel{p}{\rightarrow} \sigma^2>0$ as $n\to\infty$ if and only if, for each $\epsilon>0$,
\begin{equation}\label{eqn:cond}
\lim_{m\to \infty} \sup_{n \geq m} P\left(\frac{1}{N_n} \sum_{i=1}^{N_n} T_{b,i}^2 I(|T_{b,i}|>m) >\epsilon\right) =0.
\end{equation}
(ii) Additionally, (\ref{eqn:cond}) holds whenever $\{Y_b^2 : b\geq 1\}$  is uniformly integrable, where $Y_b$ denotes a
random variable with distribution $D_{n,b}$, $n \geq 1$, from (\ref{eqn:Dn}) (i.e., $P(Y_b \leq x) = D_{n,b}(x)$, $x \in\mathbb{R}$).  If  (\ref{eqn:2}) and $T_n = o_p(\tau_n/\tau_b)$ hold, uniform integrability of $\{Y_b^2 : b\geq 1\}$  is equivalent to $ \int x^2 d D_{n,b}(x)= N_n^{-1}\sum_{i=1}^{N_n} \E T_{b,i}^2 \rightarrow  \sigma^2$ as $n\to \infty$.\\
(iii) (\ref{eqn:cond}) also holds whenever $\{X_t\}$ is stationary and $\{T_b^2 : b\geq 1\}$  is uniformly integrable.
\end{theoremnew1}
\noindent {\sc Remark 3}: As $T_n$ is typically tight, the assumption $T_n = o_p(\tau_n/\tau_b)$ is often satisfied by a
standard condition on block length: $b\rightarrow \infty$ with  $b/n+\tau_b/\tau_n \to 0$. Block conditions are not,
in fact, used or required in statements of Theorems~\ref{theorem1}-\ref{theorem3} above. However, block assumptions
are usually needed to show the original subsampling estimator $\Ln$ is convergent as in
(\ref{eqn:2}), and examples of Section~\ref{sec:app} shall impose block length conditions for this purpose.\\

Theorem~\ref{theorem3} connects convergence (\ref{eqn:2}) of subsampling distributions $\Ln$ to the convergence
of subsampling variances  $\hatsigma$ in a way   involving no further
conditions on the process or statistic beyond  mild types of uniform integrability.
   For example, with non-stationary processes $\{X_t\}$, Theorem~\ref{theorem3}(ii)   converts the problem
of probabilistic convergence $\hatsigma \stackrel{p}{\to} \sigma^2$ into a more approachable one of subsample-moment convergence $N_n^{-1}\sum_{i=1}^{N_n} \E T_{b,i}^2 \to \sigma^2$.  To frame  another implication of Theorem~\ref{theorem3}, note that many inference problems with time series
involve a stationary process $\{X_t\}$ and a statistic $T_n$ with a normal
limit (\ref{eqn:1}) such that $\{T_n^2 : n\geq 1\}$, and consequently $\{T_b^2 : b\geq 1\}$, is uniformly integrable;
see Remark~2. In such problems, it suffices to simply  establish the consistency of the
subsampling estimator $\Ln$ (\ref{eqn:2}) and then the consistency of subsampling variance
$\hatsigma$ follows with no further effort (by Theorem~\ref{theorem3}(iii)) along with the consistency of the
convolved subsampling estimator $\Cn$ (by Theorem~\ref{theorem1}). Again, with  sample
means, $\Cn$ is a block bootstrap distribution and $\hatsigma$ is a block bootstrap variance
estimator, so both will be consistent in this setting by showing that $\Ln$ is consistent.
This strategy has two advantages with the block bootstrap: showing the consistency of
$\Ln$ is often an easier prospect than considering either $\Cn$ or $\hatsigma$ directly, and the
consistency of $\Ln$ (and thereby the bootstrap) can typically be established under weak
process assumptions.

To illustrate, Section~\ref{sec:app} applies the basic   results here for establishing the convolved subsampling estimator $\Cn$, as
  well as the block bootstrap for sample means, under differing dependence structures.

\section{Applications of convolved subsampling estimation}
\label{sec:app}
Section~\ref{sec:mix} first develops consistency results for convolved subsampling estimators with
strongly mixing processes and general statistics. The remaining subsections then consider
convolved subsampling  for the particular case of the sample mean with the goal
of generalizing and extending results for the block bootstrap across various types of dependent
data, such as non-stationary mixing time processes (Section~\ref{sec:mix2}), linear time series
(Section~\ref{sec:linear}), long-range dependent processes (Section~\ref{sec:lrd}) and spatial data (Section~\ref{sec:spatial}).
%Section~\ref{sec:convmore} draws some further connections between convolved subsampling and recent
%works with generalized subsampling methods for U-statistics and spectral estimation.

Define the strong mixing coefficient of $\{X_t\}$ as $\alpha(k) = \sup_{i\in\mathbb{Z}} \{|P(A\cap B)-P(A)P(B)|: A \in \mathcal{F}^i_{-\infty},B\in \mathcal{F}_{k+i}^\infty\}$, $k \geq 1$, where $\mathcal{F}^i_{-\infty}$ and $\mathcal{F}_{k+i}^\infty$ respectively denote $\sigma$-algebras generated
by $\{X_t : t \leq i\}$ and $\{X_t: t \geq k+i\}$ (cf.~\citep{KrishnaB.Athreya2006}, ch.~16.2). Recall $\{X_t\}$
is said to be strongly mixing or $\alpha$-mixing if $\lim_{k\to\infty} \alpha(k) = 0$.

\subsection{Convolved subsampling for general statistics under mixing}
\label{sec:mix}
For mixing stationary time series, \citet{Radu1996} proved consistency of block bootstrap
estimation for $T_n = \tau_n(t_n(X_1,\ldots,X_n) - t(P))$ based on the sample mean $t_n(X_1,\ldots,X_n)=\bar{X}_n$
 with $t(P) = \E X_1$ and $\tau_n =
\sqrt{n}$.
 The assumptions made were quite weak, requiring only
\begin{description} \itemsep 0cm
\item (a1) a stationary, $\alpha$-mixing process fulfilling (\ref{eqn:1}) (i.e., $T_n \stackrel{d}{\rightarrow} N(0,\sigma^2)$)
 and block lengths
$b^{-1}+b/n\to 0$ as $n\to\infty$;
\item (a2) uniformly integrable $\{T_n^2: n \geq 1\}$.
\end{description}
From results in Section~\ref{sec:main} and the equivalence between the block bootstrap and the convolved
subsampling estimator $\Cn$ for the sample mean, a different perspective is possible
for the bootstrap findings in \citet{Radu1996}. Under only assumption (a1) above,
the subsampling estimator $\Ln$ is consistent (i.e., (\ref{eqn:2}) holds) for the asymptotically normal
distribution of $T_n =\sqrt{n}(\bar{X}_n - \E X_1)$
  (cf.~Theorem~3.2.1, \citep{PoRoWo}),
implying, by Theorem~\ref{theorem1} here, that the block bootstrap estimator $\Cn$ would be consistent
if the subsampling variance converges $\hatsigma \stackrel{p}{\rightarrow}\sigma^2$. But, if $\Ln$ is consistent for a normal
limit by (a1), assumption (a2) then guarantees that $\hatsigma \stackrel{p}{\rightarrow}\sigma^2$ holds by Theorem~\ref{theorem3}.
Furthermore, under (a2) and with $k_n =\lfloor n/b \rfloor \rightarrow \infty$  resampled blocks as in \citet{Radu1996,Radu2012}, convergence $\hatsigma \stackrel{p}{\rightarrow}\sigma^2$ becomes even necessary here by Theorem~\ref{theorem2}. Hence,
$\alpha$-mixing serves to show that the original subsampling estimator $\Ln$ is consistent; after
which, uniform integrability and stationary assure both $\hatsigma \stackrel{p}{\rightarrow}\sigma^2$ and consistency of the
block bootstrap estimator $\Cn$ by Theorems~\ref{theorem2}-\ref{theorem3}.

Under analogously weak assumptions as those of \citet{Radu1996}, Theorem~\ref{theorem4} next
provides the general consistency of convolved subsampling estimation for general statistics
arising from mixing, and possibly non-stationary, time processes. When applied to a sample
mean $t_n(X_1,\ldots,X_n) =\bar{X}_n$, so that $\Cn$ is a block bootstrap estimator, this result extends
those of \citet{Radu1996} in two ways: by allowing potential non-stationarity series and
by permitting arbitrary levels $k_n$ of convolution/block resampling (rather than the single
choice $k_n = \lfloor n/b \rfloor$). When the statistic $t_n(X_1,\ldots,X_n)$ is not a sample mean, $\Cn$ may again
no longer correspond to the block bootstrap but has interest as a block
resampling estimator representing a composite of subsampling and bootstrap; see also Section~\ref{sec:convmore}.

\begin{theoremnew1}\label{theorem4} Let $\{X_t\}$ be a (possibly non-stationary) strongly mixing sequence. Suppose
$b^{-1}+b/n + \tau_b/\tau_n \to 0$ as $n\to \infty$; $T_n = o_p(\tau_n/\tau_b)$; (\ref{eqn:cond}) holds; and that $Y_b \stackrel{d}{\rightarrow} N(0,\sigma^2)$ as $n\to \infty$, for some $\sigma^2 > 0$, where each random variable  $Y_b$, $b\equiv b_n \geq 1$, has distribution
function $D_{n,b}$ from (\ref{eqn:Dn}).
Then, as $n \to \infty$,
\begin{align*}
\sup_{x\in\mathbb{R}} \lvert \Ln(x) - \Phi(x/ \sigma) \rvert \stackrel{p}{\rightarrow} 0 \ \ \ \ \text{and} \ \ \ \ \hatsigma \stackrel{p}{\rightarrow} \sigma^2
\end{align*}
and, for any positive integer sequence $k_n$,
\begin{align*}
\sup_{x\in\mathbb{R}} \lvert \Cn(x) - \Phi(x/ \sigma) \rvert \stackrel{p}{\rightarrow} 0.
\end{align*}
Furthermore, if (\ref{eqn:1}) additionally holds (i.e., $T_n\stackrel{d}{\rightarrow} N(0,\sigma^2)$), then $\Ln$   and $\Cn$   (with any
$k_n$) are consistent for the distribution $\Jn$ of $T_n$:
\begin{align*}
\sup_{x\in\mathbb{R}} \lvert \Ln(x) - \Jn(x) \rvert \stackrel{p}{\rightarrow} 0 \ \ \ \ \text{and} \ \ \ \ \sup_{x\in\mathbb{R}} \lvert \Cn(x) - \Jn(x) \rvert \stackrel{p}{\rightarrow} 0.
\end{align*}
\end{theoremnew1}
While providing a broad result on the validity of convolved subsampling estimation for
mixing processes, Theorem~\ref{theorem4} also expands the general subsampling results of \citet{PoRoWo}~(ch.~4.2), which focused on $\Ln$ for mixing series, to further include consistency of the
subsampling variance $\hatsigma$. That is, when dropping (\ref{eqn:cond}), the remaining Theorem~\ref{theorem4} assumptions
match those of Theorem~3.2.1-4.2.1 of \citet{PoRoWo} for the consistency of $\Ln$ to
a normal limit; including (\ref{eqn:cond}) in Theorem~\ref{theorem4} is then necessary for $\hatsigma \stackrel{p}{\rightarrow} \sigma^2$ by Theorem~\ref{theorem3}
and assures convergence of $\Cn$ by Theorem~\ref{theorem1}.

If the process $\{X_t\}$ is actually stationary, we immediately obtain the following result.

\begin{corollary1}
 \label{cor2} Let $\{X_t\}$  be a stationary, strongly mixing sequence. Suppose also $b^{-1}+b/n+\tau_b/\tau_n\to 0$ as $n \to \infty$; that (\ref{eqn:1}) holds; and that (\ref{eqn:cond}) holds (e.g., uniform integrability of
$\{T_n^2: n \geq 1\}$ suffices). Then, as $n \to\infty$, the convergence results of Theorem~\ref{theorem4} hold.
\end{corollary1}
Section~\ref{sec:mix2} next provides some further refinements with mixing processes in the sample
mean case, where $\Cn$ matches the block bootstrap.

\subsection{Block bootstrap for mixing non-stationary time processes}
\label{sec:mix2}
Consider a strongly mixing, potentially non-stationary   sequence $\{X_t\}$ having a common mean
parameter $\E X_t = \mu\in\mathbb{R}$,   which is  estimated by the sample mean $\bar{X}_n$.
 In this setting and under conditions where $ T_n \equiv \sqrt{n}(\bar{X}_n-\mu)$  has a
normal limit (\ref{eqn:1}), \citet{Fitz}  established the consistency of the block bootstrap for estimating
the distribution of $T_n$.  The result, however, required the existence of a $(4 + \delta)$-moment (i.e., $\sup_t \E |X_t|^{4+\delta}<\infty$ for some $\delta > 0$)  along with stringent mixing conditions and
restrictions on the block length $b = o(n^{1/2})$.
\citet{PoRoWo}~(example 4.4.1) showed
that the subsampling estimator $\Ln$ is consistent under  weaker conditions, including
only a $(2 + \delta)$-moment.
For the block bootstrap with $k_n = \lfloor n/b\rfloor$ resampled blocks, \citet{PoRoWo} also proved bootstrap consistency by applying convolved subsampling in this problem, using a weaker block assumption  $b = o(n)$ than \citet{Fitz}  but otherwise
with same remaining strong assumptions about the process. However,
\citep{PoRoWo}~(remark 4.4.4) conjectured that the block bootstrap might be established
under non-stationary
using the same weak moment/mixing conditions as the subsampling estimator
$\Ln$, just as in the case of stationary mixing processes (cf.~\citep{Radu1996}). We confirm
this by the following Theorem~\ref{theorem5}.

\begin{theoremnew1} \label{theorem5}
Let $\{X_t\}$ be a sequence of (not necessarily stationary) strongly mixing random
variables with common mean $\mu$. For some $\delta > 0$, suppose that $\sup_t \E |X_t|^{2+\delta}<\infty$   and
$\sum_{k=1}^\infty \alpha(k)^{\delta/(2+\delta)}<\infty$.  Assume also that, for some $\sigma^2>0$,
\[
\lim_{n\to \infty} \sup_{i \geq 1} \left|\var\left(n^{-1/2} \sum_{t=i}^{i+n-1} X_t \right)-\sigma^2 \right|=0.
\]
Then, as $n\to \infty$, $T_n=\sqrt{n}(\bar{X}_n-\mu) \stackrel{d}{\rightarrow} N(0,\sigma^2)$ (i.e., (\ref{eqn:1}) holds).  Additionally, if $b^{-1}+b/n\rightarrow 0$ as $n\to \infty$, then
\[
\sup_{x \in \mathbb{R}}  \lvert \Ln(x) - \Phi(x/ \sigma) \rvert \stackrel{p}{\rightarrow} 0 \ \ \ \ \text{and} \ \ \ \ \hatsigma \stackrel{p}{\rightarrow} \sigma^2
\]
and, for any positive integer sequence $k_n$,
  \[
\sup_{x \in \mathbb{R}} |\Cn(x) - \Phi(x/ \sigma)| \stackrel{p}{\rightarrow} 0.
\]
\end{theoremnew1}
Hence, with any number $k_n$ of concatenated blocks, the block bootstrap estimator $\Cn$ is
valid for the distribution of the sample mean under  mild assumptions for mixing, and
possibly non-stationary, processes. Note that the assumptions of Theorem~\ref{theorem5} resemble those
essentially needed for a central limit theorem (CLT) itself (cf.~Theorem~16.3.5, \citep{KrishnaB.Athreya2006}). In particular, the assumptions also match those commonly used in the stationary
case for establishing the block bootstrap; see Section 3.2 of \citet{Lah2}. With the same
moment condition as \citet{PoRoWo}~(Theorem~4.4.1), Theorem~\ref{theorem5} additionally shows that
the original subsampling estimator $\Ln$ is consistent under non-stationarity with even
weaker mixing assumptions than considered previously
$\sum_{k=1}^\infty
(k +1)^2 \alpha(k)^{\delta/(8+\delta)}<\infty $.  The
central message of Theorem~\ref{theorem5}, however, is that the convolved subsampling approach allows
the block bootstrap estimator $\Cn$ for the sample mean to be established under weak
conditions similarly to $\Ln$.

Next consider the block bootstrap in another important example of non-stationarity,
involving certain periodically correlated time series.  Here
the mean function $\mu(t) \equiv \E X_t$ is not constant, as in Theorem~\ref{theorem5}, but rather an almost periodic function.
A real-valued function $f$ is \emph{almost periodic} if, for every $\epsilon>0$, there is an $n(\epsilon)\in\mathbb{N}$ such that in every interval $I_{n(\epsilon)}$ of length $n(\epsilon)$ or greater, there is an integer $p \in I_{n(\epsilon)}$  such that
\begin{align*}
\sup_{t \in \mathbb Z} \lvert f(t+p) -f(t) \rvert < \epsilon;
\end{align*}
see \citep{Cor}. For such functions the limit $M(f) \equiv \lim_{n \to \infty} n^{-1} \sum_{i=s}^{s+n-1}f(i)$ exists and does not depend on $s$. Moreover, if the set $\Lambda =\{ \lambda \in [0,2\pi): M( g_\lambda) \not = 0 \}$ is finite for $g_\lambda(t) \equiv f(t) e^{- \imath \lambda t}$, $t\in\mathbb{Z}$ ($\imath=\sqrt{-1}$), then
\begin{align}
 \left|\frac 1 n \sum_{i=s}^{s+n-1}(f(i) - M(f) )  \right|\leq \frac C n \label{InequalAlmostPeriodic}
\end{align}
holds for some $C>0$ not depending on $n$ or $s$ by \citet{CaHoHuLe}. Hence,
$M(f)$ represents the mean value of an almost periodic function $f$. A time series is called \emph{almost periodically correlated} (APC) if its mean and autocovariance functions are almost
periodic, that is, for  every fixed $\tau \in \mathbb Z$,
\begin{align*}
\mu(t) = \E X_t \ \ \ \ \text{and} \ \ \ \ \rho_{\tau}(t) = \E X_t X_{t + \tau}
\end{align*}
are almost periodic as functions of $t$; see \citep{Hur}. For an ACP series $\{X_t\}$, a parameter of interest is then $t(P)\equiv M( \mu) = \lim_{n \to \infty} n^{-1} \sum_{i=s}^{s+n-1} \mu(i)$
 as a summary of the
process mean structure, which is estimated by  $\bar{X}_n$.
 \citet{Syn} showed that the block bootstrap consistently estimates the sampling
distribution of $T_n = n^{1/2}( \bar{X}_n-M(\mu))$
  under appropriate conditions.  By applying the
convolved subsampling technique, we may extend the   bootstrap results of \citet{Syn} (Corollary 3.2) by  substantially weakening the  assumptions made
there about $(4 + \delta)$-moments and
$\sum_{k=1}^\infty k\alpha(k)^{\delta/(4+\delta)}<\infty$.

\begin{corollary1}
\label{cor3}
Let $\{X_t\}$ be an APC sequence of strongly mixing random variables such that $\sup_{t } \E |X_t|^{2+\delta}<\infty$ and   $\sum_{k=1}^{\infty}   \alpha(k)^{\delta/(\delta+2)}  < \infty$
for some $\delta>0$, and suppose
the set $\Lambda =\{ \lambda \in [0,2\pi): M(g_\lambda) \not = 0 \}$ is finite for $g_\lambda(t) \equiv \mu(t)e^{-\imath\lambda t}$, $t\in\mathbb{Z}$, with $\mu(t)=\E X_t$. Then,
all conclusions of Theorem~\ref{theorem5} hold for $T_n = n^{1/2}( \bar{X}_n-M(\mu))$  as $n\to \infty$.
\end{corollary1}

\subsection{Block bootstrap for linear time processes}
\label{sec:linear} Based on a sample $X_1,\ldots,X_n$, next consider inference about the mean $\E X_t = \mu\in\mathbb{R}$   of a stationary  time
process $\{X_t\}$ prescribed as
\begin{equation}\label{eqn:linear}
X_t = \mu  + \sum_{j\in\mathbb{Z}} a_j \varepsilon_{t-j},\quad t\in\mathbb{Z},
\end{equation}
in terms of iid variables $\{\varepsilon_t\}$ with mean zero and finite variance $\E \varepsilon_t^2 \in(0,\infty)$  and a real-valued
sequence $\{a_j\}$ of constants where $\sum_{j\in\mathbb{Z}}a_j^2<\infty$.
  The linear series $\{X_t\}$
need not be mixing and, depending on constants $\{a_j\}$, can potentially exhibit either weak
or strong forms of time dependence. Considering the sample mean $\bar{X}_n$
as an estimator of
the process mean $\mu$, suppose that
\begin{equation}\label{eqn:var}
\lim_{n\to \infty} n^{\alpha} \var(\bar{X}_n) = \sigma^2
\end{equation}
for some $\sigma^2 > 0$ and exponent $\alpha\in  (0, 1]$ depending on the process $\{X_t\}$. When
$\alpha=1$, the sample mean's variance    decays at a rate $O(n^{-1})$ with the sample size, as
typical for weakly, or short-range, dependent processes. However, when $\alpha \in (0,1)$,
the sample mean has a variance with comparatively slower decay $O(n^{-\alpha})$, which may be
associated with processes exhibiting strong, or long-range, forms of dependence. Long-range
dependent processes are commonly characterized by slowly decaying covariances involving
a long-memory exponent $\alpha\in (0, 1)$, which  results in less precision (\ref{eqn:var}) for a sample
mean compared to the weak dependence case \citep{Beran94}. Classes of strongly dependent
 processes  that satisfy (\ref{eqn:linear})-(\ref{eqn:var})  include
fractional Gaussian models \citep{MVN68} and fractional autoregressive
integrated moving averages \citep{GJ80, Hosking81}.

Based on (\ref{eqn:var}), define $T_n \equiv n^{\alpha/2}(\bar{X}_n-\mu)$   in terms of scaling $\tau_n\equiv n^{\alpha/2}$. In this setting,
the convolved subsampling $\Cn$ once again corresponds to the block bootstrap estimator
based on $k_n$ resampled blocks, but there is a wrinkle to note. Recalling from (\ref{eqn:bb}) that the
bootstrap sample mean $\bar{X}^*_{n1}$
  is created from a
bootstrap sample of length $n_1 = k_n b$,
the
bootstrap analog of $T_n$ here is given by
\begin{equation}
\label{eqn:bb2}
T_n^* \equiv b^{(1-\alpha)/2} (n_1)^{\alpha/2} (\bar{X}^*_{n_1} - \E_* \bar{X}^*_{n_1})
\end{equation}
rather than $T_n^* = (n_1)^{\alpha/2} (\bar{X}^*_{n_1} - \E_* \bar{X}^*_{n_1})$. While intuitive, the latter definition is incorrect
and produces a degenerate bootstrap limit,   shown by \citet{Lah}. Instead, scalar adjustment
$b^{(1-\alpha)/2}$  is required in the bootstrap version (\ref{eqn:bb2}) of $T_n$, which disappears under weak dependence
$\alpha = 1$ whereby bootstrap versions in (\ref{eqn:bb}) and (\ref{eqn:bb2}) then match. Interestingly,
convolved subsampling estimator $\Cn$ automatically corresponds to the correct bootstrap
rendition $T_n^*$ in (\ref{eqn:bb2}) under both weak $\alpha = 1$ and strong $\alpha \in (0, 1)$ dependence.

Considering the sample mean from stationary linear processes (\ref{eqn:linear}) ranging over
short- or long-range dependence, \citet{KiNo}  showed the consistency of
the block bootstrap distribution $\Cn$ (when $k_n = \lfloor n/b \rfloor$)  and bootstrap variance $\hatsigma$.
Via convolved subsampling, we may generalize their results.  For linear processes $\{X_t\}$ satisfying (\ref{eqn:linear})-(\ref{eqn:var}),  the sample mean $T_n\equiv n^{\alpha/2}(\bar{X}_n-\mu)$
has a normal limit (\ref{eqn:1}) (cf.~\citep{Dav70}), and \citet{LaNo}  showed the
consistency of the subsampling estimator $\Ln$ (i.e., (\ref{eqn:2}) holds)
under mild assumptions. Hence, by primitively assuming (\ref{eqn:1})-(\ref{eqn:2}) to hold, Corollary~\ref{cor4} next
extends the block bootstrap to general stationary processes with sample means satisfying a
variance condition (\ref{eqn:var}), which includes results of  \citet{KiNo} for   linear processes as a special case.
\begin{corollary1} \label{cor4} Let $\{X_t\}$  be a stationary process with mean $\mu \in \mathbb{R}$ satisfying (\ref{eqn:var}) for some
$\alpha \in  (0, 1]$, and suppose that (\ref{eqn:1})-(\ref{eqn:2})  hold for $T_n\equiv n^{\alpha/2}(\bar{X}_n-\mu)$. Then, as $n\to \infty$,
\[
\hatsigma \stackrel{p}{\rightarrow} \sigma^2 \ \ \ \ \text{and} \ \ \ \ \sup_{x \in \mathbb{R}} |\Cn(x) - \Phi(x/ \sigma)| \stackrel{p}{\rightarrow} 0
\]
 for any positive integer sequence $k_n$.
\end{corollary1}
Corollary~\ref{cor4} is an application of Theorems~\ref{theorem1} and~\ref{theorem3} for stationary processes which may
not be mixing. Our exposition has assumed the exponent $\alpha \in   (0, 1]$ to be known. Upon
replacing $\alpha$ with an estimator $\hat{\alpha} \equiv \hat{\alpha}(X_1,\ldots,X_n)$  where $|\hat{\alpha}-\alpha| \log n \stackrel{p}{\rightarrow} 0$, the conclusions
of Corollary~\ref{cor4} still hold; see Remark 3 of \citep{KiNo} for further details.
\subsection{Block bootstrap under long-range dependence}\label{sec:lrd}
This section briefly mentions the block bootstrap with additional types of long-memory
sequences. Beyond linear processes, the sample mean of a long-range dependent sequence
may converge to a non-normal limit, such as the case for certain subordinated Gaussian
processes considered by \citet{Taq,Taq3} and \citet{DoMa} (e.g., $X_t =
G(Z_t)$ as a function $G$ of a long-range dependent Gaussian series $\{Z_t\}$). For such time series, \citet{Lah} proved that the block bootstrap sample mean always has a normal limit, so
that the block bootstrap fails if the original sample mean is asymptotically non-normal.
This result is in concordance with our Theorem~\ref{theorem2}(ii).

\citet{HoWeWuZh} considered subsampling estimation for a wider class of long-range
dependence series that includes both subordinated Gaussian processes as well types of
linear processes (\ref{eqn:linear}). Namely, sequences $X_t = K(Z_t)$, $t\in\mathbb{Z}$, formed by a measurable
transformation $K$ of a long-range dependent linear process
\begin{align*}
Z_t =  \varepsilon_t + \sum_{j=1}^{\infty} j^{-\beta} L(j)  \varepsilon_{t-j}, \quad t\in\mathbb{Z},
\end{align*}
defined with iid mean zero, finite variance innovations $\{\varepsilon_t\}$, an index parameter $1/2 < \beta < 1$,
and slowly varying function $L(\cdot)$. They distinguish two cases, depending on $\beta$ and the so-called
power rank $p \geq 1$ of $K$. In the first case (i.e., $p(2\beta -1) > 1$), the transformation
$K$  diminishes long-range dependence,  and the sample mean converges to
a normal limit. In the second case (i.e., $p(2\beta- 1) < 1$), the transformed process
$X_t = K(Z_t)$ remains strongly dependent and  the sample mean has a normal limit only
when $p = 1$.

Assuming a constant function $L(\cdot) = C$ in the above formulation, the variance
of a sample mean satisfies (\ref{eqn:var}) (i.e., $\lim_{n\to \infty} n^\alpha \var(\bar{X}_n)=\sigma^2>0$)
  with a long-memory
exponent $\alpha \equiv \min \{1, p(2\beta-1) \} \in (0, 1]$ that changes between cases of weak $\alpha = 1$ or strong
$\alpha = p(2\beta-1) \in (0, 1)$ dependence (cf.~Lemma 1, \citep{HoWeWuZh}). For the sample
mean, \citet{HoWeWuZh} established consistency of several subsampling estimators as well
as convergence of  $\hatsigma$. Thus, by slightly re-casting results of \citep{HoWeWuZh} and applying our Corollary~\ref{cor4}, we may show the validity of the block bootstrap
$\Cn$ for estimating the distribution of $T_n \equiv n^{\alpha/2}( \bar{X}_n - \mu)$, $\mu= \E X_t$,
  for transformed linear
processes exhibiting either short- or long-range dependence. To the best of our knowledge,
 the bootstrap has not yet been investigated   for such processes.
 \begin{corollary1} \label{cor5} For $X_t=K(Z_t)$, $t\in\mathbb{Z}$, as above,
  suppose (\ref{eqn:var}) holds for $\alpha= \min\{1, p(2\beta- 1)\} \in   (0, 1]$ along with conditions
of Theorem~1 in \citep{HoWeWuZh} (involving a block $b \propto n^a$ for some $a \in (0, 1)$) with
either $p(2\beta- 1) > 1$, or $p = 1$ and $(2\beta- 1) < 1$. Then, for $T_n = n^{\alpha/2}( \bar{X}_n - \mu)$ as $n\to \infty$,
  both  (\ref{eqn:1})-(\ref{eqn:2})  hold and
\[
\hatsigma \stackrel{p}{\rightarrow} \sigma^2 \ \ \ \ \text{and} \ \ \ \ \sup_{x \in \mathbb{R}} |\Cn(x) - \Phi(x/ \sigma)| \stackrel{p}{\rightarrow} 0
\]
 for any positive integer sequence $k_n$.
\end{corollary1}
As with Gaussian subordinated processes \citep{Lah},   consistency of the
block bootstrap or convolved estimator $\Cn$ holds only in cases where the sample
 mean  follows a CLT. For statistics other  than the
sample mean,   \citet{BeWe} proved the consistency of the subsampling estimator $\Ln$ with transformations of Gaussian processes, provided
that the original statistic converges under mild assumptions. When this limit is normal, consistency of the convolved
estimator will follow by our Theorem~\ref{theorem2} by showing convergence of $\hatsigma$ (which, as \citet{BeWe} consider only stationary processes, can hold by Theorem~\ref{theorem3} and uniform integrability of the statistic).

\subsection{Spatial data}
\label{sec:spatial}
 While   convolved subsampling results have been presented
 for processes $\{X_t\}$ indexed by time $t$ to ease the exposition,  Theorems~\ref{theorem1}-\ref{theorem3} also  apply
to more general processes, including spatial random fields.    We illustrate this with spatial data on a grid, for which  various authors have considered
block bootstrap and
subsampling methods;  see \citet{PoRoWo} (ch.~5) and \citet{Lah2} (ch.~12) and references therein. For approximating spatial sample means, we
examine convolved subsampling and the block bootstrap  under particularly weak conditions.

To describe the spatial sampling, suppose $\{X(\bm{s}):\bm{s}\in\mathbb{Z}^d\}$ denotes a stationary random field which is observed at $n$ locations $\{\bm{s}_1\ldots,\bm{s}_n\} = R_n \cap \mathbb{Z}^d$
within a spatial sampling region $R_n \subset \mathbb{R}^d$.  Here $d \geq 1$ denotes the dimension of sampling
and  spatial variables have regular locations denoted by the integer grid $\mathbb{Z}^d$.  A common formulation of the sampling region $R_n = \lambda_n R_0$  involves inflating a template set $R_0 \subset (-1/2,1/2]^d$ by a sequence of positive scaling factors $\lambda_n\rightarrow \infty$
as $n\to \infty$ (cf.~\citep{Cre}), where $R_0$
  contains the origin and satisfies a mild boundary condition
to avoid pathological geometries  (cf.~Condition~B, ch.~12.2, \citep{Lah2}).

  Based on the
sample mean $t_n(X_1,\ldots,X_n)=\bar{X}_n$ of the data $X(\bm{s}_i)\equiv X_i$, $i=1,\ldots,n$ and   the process mean $t(P)\equiv \E X(\bm{s})=\mu$, we   approximate
the distribution of $T_n = \sqrt{n}(\bar{X}_n-\mu)$ through spatial subsamples, described next.
Using a  template $\mathcal{D}_0\subset (-1/2,1/2]^d$ (similar to $R_0$) and a positive sequence $b\equiv b_n$ of subsample scaling such that $b^{-1}+ b/\lambda_n\rightarrow 0$ as $n\to \infty$, define a prototype spatial subregion  $\mathcal{D}_b \equiv b $ and an  index set $\mathcal{I}_n \equiv \{\bm{i}\in \mathbb{Z}^d : \bm{i}+ \mathcal{D}_b\subset R_n\}$ of all integer translates of   $\mathcal{D}_b$ lying inside $R_n$.  Each subsampling region $\bm{i}+\mathcal{D}_b$, $\bm{i}\in\mathbb{Z}^d$,  provides an analog statistic $t_{n,b,\bm{i}} \equiv \sqrt{n_b} ( \bar{X}_{b,\bm{i}}-\bar{X}_n)$ from the average $\bar{X}_{b,\bm{i}} \equiv \sum_{\bm{s}\in \bm{i}+\mathcal{D}_b}X(\bm{s})/n_b$ of the $n_b=|\mathcal{D}_b\cap \mathbb{Z}^d|$ observations per subregion.
Re-writing $\{t_{n,b,\bm{i}}: \bm{i}\in\mathcal{I}_n\}$ as $\{t_{n,b,i}\}_{i=1}^{N_n}$ for notational simplicity, with $N_n=|\mathcal{I}_n|$ denoting the number of subsamples, the subsampling estimator $\Ln$
 is then computed as in (\ref{eqn:sub}), inducing also a convolved estimator $\Cn$.
Above, the    template   $\mathcal{D}_0$  can have an arbitrary shape, with the choice
$\mathcal{D}_0=R_0$ producing subregions as scaled-down copies of $R_n$.  However, in selecting a cube $\mathcal{D}_0=(-1/2,1/2]^d$, the  convolved estimator $\Cn$ then matches a standard spatial block bootstrap
for the sample mean based on $k_n$ resampled blocks; see  \citep{Lah2} (ch.~12.3) and \citep{nordlahfrid}.
%In which case,  $k_n$
%is typically taken as the number of disjoint cubes
%$b (-1/2,1/2]^d$ needed to tile $R_n$, so that $k_n \approx n/n_b\to \infty$ holds similarly to time series (where $k_n \approx n/b $ for $n_b=b$).

Under appropriate  assumptions for the stationary  random field,  the  spatial sample mean $T_n =\sqrt{n}(\bar{X}_n-\mu)$ has a normal limit; see \citep{Lah03} for references and a general development.  We establish convolved  subsampling under mixing conditions from  \citet{Lah03} (sec.~4.2) which are almost optimal, or  minimal,  for a spatial CLT.
For $a>0$ and $c \geq 1$, define the mixing coefficient $\alpha(a;c)$ of   $\{X(\bm{s}):\bm{s}\in \mathbb{Z}^d\}$ as
\[
 \sup\{|P(A_1 \cap A_2)-P(A_1)P(A_2)|: A_i\in\mathcal{F}(S_i), S_i\in V(c), i=1,2; \mathrm{dist}(S_1,S_2) \geq a\},
\]
where $V(c)$ is the collection of all sets in $\mathbb{R}^d$, with a volume of $c$ or less, that can be expressed as unions of up to $\lceil c \rceil$ many cubes; for a set $S \subset \mathbb{R}^d$, $\mathcal{F}(S)$ is the $\sigma$-algebra generated by $\{X(\bm{s}):\bm{s}\in S\cap \mathbb{Z}^d\}$; and
$\mathrm{dist}(S_1,S_2) = \inf\{|\bm{s}_1-\bm{s}_2|_\infty:\bm{s}_1\in S_1, \bm{s}_2\in S_2 \}$
denotes the distance between sets $S_1$ and $S_2$ with
$|\bm{s}|_\infty = \max_{1 \leq i \leq d}|s_i|$ for $\bm{s}=(s_1,\ldots,s_d)^\prime \in \mathbb{R}^d$.

\begin{theoremnew1}
\label{theorem6} Suppose that the random field $\{X(\bm{s}):\bm{s}\in\mathbb{Z}^d\}$ is stationary with
 $\E |X(\bm{0})|^{2+\delta}<\infty$ and with mixing coefficient satisfying
 \[
 \alpha(a;c) \leq C a^{-\tau_1} c^{\tau_2},\quad a,c \geq 1,
 \]
 for some $\delta>0$, $C>0$, $\tau_1> d(2+\delta)/\delta$ and $0 \leq \tau_2 \leq (\tau_1-d)/(4d)$, and   assume $\sigma^2 \equiv \sum_{\bm{k}\in \mathbb{Z}^d} \mathrm{Cov}(X(\bm{0}),X(\bm{k})) >0$. Then, as $n\to \infty$,  $T_n \equiv \sqrt{n}(\bar{X}_n-\mu) \stackrel{d}{\rightarrow}N(0,\sigma^2)$  (i.e., (2) holds).  Additionally, if $b^{-1}+ b/\lambda_n\rightarrow 0$ as $n\to \infty$, then
 \begin{align*}
\sup_{x\in\mathbb{R}}\lvert \Ln(x) - \Phi(x/ \sigma) \rvert \stackrel{p}{\rightarrow} 0 \ \ \ \ \text{and} \ \ \ \   \hatsigma\stackrel{p}{\rightarrow} \sigma^2
\end{align*}
and, for any positive integer sequence $k_n$,
\begin{align*}
\sup_{x\in\mathbb{R}} \lvert \Cn(x) - \Phi(x/ \sigma) \rvert \stackrel{p}{\rightarrow} 0.
\end{align*}
\end{theoremnew1}

Theorem~\ref{theorem6} consequently demonstrates the spatial block bootstrap for the sample mean (i.e., template $\mathcal{D}_0=(-1/2,1/2]^d$) under weaker conditions than considered previously, such as $ \E |X(\bm{0})|^{6+\delta}<\infty$ and $\tau_1 > 5d(6+\delta)/\delta$ (cf.~Theorem~12.1, \citet{Lah2}).   While such earlier assumptions  were made to ease proofs, convolved subsampling    also allows for simple proofs under  mild spatial assumptions.  The mixing assumptions in
Theorem~\ref{theorem6} are also comparable to those used by \citet{PoRoWo} (Theorem~5.3.1) for showing the consistency of the subsampling distribution $\Ln$
with rectangular sampling regions $R_n$.

%Several authors have also considered block resampling methods for irregularly located spatial data,
%such as \citet{PoPaRo}, \citet{polsher} and \citet{lahzhu06}. Convolved subsampling may also be extended
%to such spatial data, though the spatial asymptotics involved require a more complicated development.

\section{Convolved subsampling   in other contexts}
\label{sec:convmore}
We briefly outline relationships between convolved subsampling and some recent works
with the block resampling for statistics which fall outside of the sample mean cases in Sections~\ref{sec:mix2}-\ref{sec:spatial}.
These works, presented next in Section~\ref{sec:ustat} for U-statistics and Section~\ref{sec:cov} for spectral estimators,
involve mixing time processes, and so can be connected to the general convolved subsampling result (Theorem~\ref{theorem4}) of Section~\ref{sec:mix} for arbitrary statistics.

\subsection{U-statistics}

\label{sec:ustat}
U-statistics are a class of nonlinear functionals for prescribing   statistics, such as the sample variance.  Suppose that $X_1,\ldots,X_n$ arise from a stationary process and, based on a symmetric kernel $h:\mathbb{R}^2\to \mathbb{R}$,
define a (bivariate) U-statistic as
\[
t_n\equiv t_n(X_1,\ldots,X_n) = \frac{2}{n(n-1)} \sum_{1 \leq i<j\leq n} h(X_i,X_j),
\]
 which estimates a target parameter $t(P)\equiv \int h(x,y) d G(x) dG(y)$, where $G$ denotes the marginal distribution of $X_t$.  Consider the problem of estimating the distribution of
 $T_n\equiv \sqrt{n}(t_n-t(P)) $, with scaling $\tau_n=\sqrt{n}$, under weak time dependence.
  The subsampling distribution   $\Ln$ is defined  by computing the U-statistic
  $t_{n,b,i} =  [b(b-1)]^{-1}2\sum_{i \leq j_1<j_2\leq i+b-1} h(X_{j_1},X_{j_2}) $ on each length $b$ subsample $\{(X_i,\ldots,X_{i+b-1})\}_{i=1}^{N_n\equiv n-b+1} $  in (\ref{eqn:sub}).     In contrast, block bootstrap versions of U-statistics have a formulation similar to (\ref{eqn:bb}); see \citet{DeWe}, \citet{ShWe} and \citet{Leu}. That is,  a bootstrap sample  $X_1^{\ast}, \dots , X_{n_1}^{\ast}$, $n_1=k_n b$, is  generated by
  resampling $k_n$ blocks of length $b$ (typically $k_n =\lfloor n/b \rfloor$) and then the U-statistic $t_{n_1}^*\equiv t_{n_1}(X_1^*,\ldots,X_{n_1}^*)$ is calculated from the complete bootstrap sample to create a bootstrap rendition $T_n^* = \sqrt{n_1} (t_{n_1}^*- \E_* t_{n_1}^*)$ of $T_n$.
In this setting, the bootstrap distribution $T_n^*$ would not generally correspond that of a $k_n$-fold convolution $\Cn$ of the subsampling distribution $\Ln$, as occurred in the sample mean case (Section~\ref{sec:conv}).

However,  \citet{ShTeWe2} recently considered an alternative block resampling estimator for U-statistics,
which matches the convolved subsampling estimator $\Cn$ here based on the subsampling estimator $\Ln$ for $T_n$ described above.
 Note that, for stationary mixing data,    \citet{DeWe} (Theorem~1.8-Lemma~3.6) provide a CLT for the relevant U-statistic:
 $T_n \stackrel{d}{\to} N(0,\sigma^2)$ and $\E T_n^2 \to \sigma^2$ as $n\to \infty$ where $\sigma^2 \equiv 4 \sum_{k=-\infty}^\infty \mathrm{Cov}(h_1(X_0),h_1(X_k))$ for $h_1(x) = \int h(x,y) d G(y) $.
Under mixing conditions and with  $k_n=\lfloor n/b \rfloor\to \infty$, \citet{ShTeWe2} established that $\Cn$
captures this limiting normal distribution of $T_n$ and also
 showed the consistency of the variance $\hatsigma$ of $\Cn$.  The argument there involved decomposing the bootstrap U-statistic $T_n^*$
into a linear part, coinciding with a sample mean  from the usual block bootstrap, and degenerate part shown to be negligible.
 However, the general convolution result  in Theorem~\ref{theorem4} for mixing processes provides an alternative, and much simpler, approach.  From $T_n \stackrel{d}{\to} N(0,\sigma^2)$ and $\E T_n^2 \to \sigma^2$,
  all of the conditions of Theorem~\ref{theorem4} automatically hold,
 proving that $\Cn$ is consistent for the distribution of $T_n$ for any convolution level $k_n$ and also that $\hatsigma \stackrel{p}{\to} \sigma^2$. This approach also weakens the block assumptions used by \citep{ShTeWe2} (i.e., $b=O(n^\epsilon)$ for some $\epsilon\in(0,1)$) to
 $b^{-1}+b/n\to 0$ under Theorem~\ref{theorem4}.

As a side note,
simulations in \citep{ShTeWe2} also indicate that convolved subsampling $\Cn$ exhibits better coverage accuracy than traditional subsampling $\Ln$ for U-statistics,
and has performance comparable to the standard block bootstrap.   This finding suggests consideration of convolved subsampling estimators even when these do not match the  block bootstrap.\\

\noindent {\sc Remark 4}:  Convolved subsampling can reduce skewness compared
to distributional estimates from subsampling; see   \citep{PoRoWo} (sec.~10.2).  This aspect may partially explain the better U-statistic coverage mentioned above.
    For example, in approximating sample means, the distribution of $T_n=\sqrt{n}(\bar{X}_n-\mu)$
often has approximate skewness $ \gamma/\sqrt{n}$ for some constant $\gamma$.  In this case, the corresponding subsampling estimator $\Ln$
has a larger approximate skewness $ \gamma/\sqrt{b}$, while a more fully convolved estimator (bootstrap) $\Cn$ with $k_n \approx n/b$
has skewness approximately $ \gamma /\sqrt{n}$, where a better matching skewness can  impact   higher-order accuracy.

\subsection{Spectral estimators for non-stationary time series}
\label{sec:cov}
As described in Section~\ref{sec:mix2}, almost periodically correlated (APC) time series $\{X_t\}$ are an important example of non-stationary sequences. Beyond the mean function, inference about the correlation structure is also of interest.  Based on a sample $X_1,\ldots,X_n$, a symmetric kernel $w(\cdot)$ and a bandwidth choice $L_n$,
  \citet{Len,Len2} considered kernel  estimators
\[
t_n(X_1,\ldots,X_n) \equiv \frac{1}{2\pi n} \sum_{t=1}^n \sum_{s=1}^n  \frac{1}{L_n} w\left(\frac{t-s}{L_n}\right) X_t X_s e^{-\imath \upsilon t}e^{ \imath \omega s}
\]
for an extended spectral density $t(P)\equiv t(P)(\upsilon,\omega)$, $(\upsilon,\omega)\in(0,2\pi]^2$, used to represent the almost periodic covariance function
$c_\tau(t)\equiv\mathrm{Cov}(X_t, X_{t+\tau})$, $t\in\mathbb{Z}$, for a given $\tau\in\mathbb{Z}$; see \citep{Len,Len2} for details.

For $T_n \equiv \tau_n (t_n - t(P))$ with scaling $\tau_n =\sqrt{n/L_n}$,
\citet{Len} (Theorems~3.1-3.2) proved a CLT $T_n \stackrel{d}{\to}N(0,\sigma^2)$ and moment convergence $  \E T_n^2 \to \sigma^2$ with mixing APC series, which was extended in \citet{Len2} to multivariate data.   Due to the complicated form of $\sigma^2$, a subsampling estimator $\Ln$
for the distribution of $T_n$
may computed as
in (\ref{eqn:sub}) with analog statistics $t_{n,b,i}$ and  scaling $\tau_b =\sqrt{b/L_b}$ defined from subsamples $\{(X_i,\ldots,X_{i+b-1})\}_{i=1}^{N_n\equiv n-b+1}$.  \citet{Len} proved the consistency of the estimator $\Ln$, while
   \citet{Len2} proposed a generalized resampling method which essentially corresponds
     a convolved subsampling estimator $\Cn$ induced from $\Ln$  (though  \citep{Len2} also considered
       $\Ln$ defined with possibly non-uniform weights on $\{t_{n,b,i}\}_{i=1}^{N_n}$).  In particular, \citet{Len2}  established the consistency of
     $\Cn$ through bootstrap arguments requiring much stronger mixing and moment assumptions than needed for the convergence
       $T_n \stackrel{d}{\to}N(0,\sigma^2)$ and  $  \E T_n^2 \to \sigma^2$.  However, the general convolved subsampling result in Theorem~\ref{theorem4}  may   alternatively be used here with mixing non-stationary ACP series.

        To apply  Theorem~\ref{theorem4} with   blocks  where $b^{-1}+b/n+\tau_b/\tau_n\to 0$ as $n\to\infty$,
       one requires that  $Y_b \stackrel{d}{\to}N(0,\sigma^2)$     and that (\ref{eqn:cond}) holds, where $Y_b$, $b\equiv b_{n} \geq 1$,
         denotes a sequence of variables with   distribution $D_{n,b}(\cdot)$ from (\ref{eqn:Dn}).
                    But, the same conditions needed for $T_n \stackrel{d}{\to}N(0,\sigma^2)$ and  $  \E T_n^2 \to \sigma^2$
         also yield  $Y_b \stackrel{d}{\to}N(0,\sigma^2)$   and $  \E Y_b^2 \to \sigma^2$ (cf.~Theorems~3.1-3.2 and 4.1, \citep{Len}).
         Furthermore, mixing and $Y_b \stackrel{d}{\to}N(0,\sigma^2)$, along with $T_n=O_p(1)$ and $\tau_n/\tau_b\to\infty$, guarantee that (\ref{eqn:2})
         holds (i.e., $\sup_{x\in\mathbb{R}}|\Ln(x)-\Phi(x/\sigma)|\stackrel{p}{\to} 0$) and that consequently
           (\ref{eqn:cond})  follows from Theorem~\ref{theorem3}(ii) by $  \E Y_b^2 \to \sigma^2$.
           That is, the same minimal conditions for a CLT with APC series  suffice for the consistency of convolved subsampling   $\Cn$ by the general result of Theorem~\ref{theorem4}.

\section{Independent data versions}
\label{sec:indep}
For completeness, we briefly mention a variation of convolved subsampling appropriate for independent data.
Recall Section~\ref{sec:mix} considered block-based
convolved subsampling with general statistics computed from strongly mixing time processes $\{X_t\}$.
Hence, results of Section~\ref{sec:mix} apply to immediately independent data, as do block bootstrap results of Section~\ref{sec:mix2} for sample means under mixing conditions.  However, with independent $X_1,\ldots,X_n$, one may consider
a different formulation of subsamples  rather than  data blocks   of   $b$ consecutive observations.    Namely, let $b_n \equiv b$ denote a   set size (rather than a block length) and define subsamples $Y_{b,1},\ldots,Y_{b,N_n}$ as the $N_n \equiv {n \choose b}$  unordered  subsets of size $b$ from $\{X_1,\ldots,X_n\}$.  The ``independent data" subsampling estimator $\Lni$ is   defined as $\Ln$ in (\ref{eqn:sub}) with   subsample statistics $t_{n,b,i} \equiv t_b(Y_{b,i}), $ $i=1,\ldots, N_n$, where   statistics $t_b(\cdot)$ are  symmetric in their arguments here.  For iid data in particular, see \citet{PoRoWo} (ch.~2) for a general treatment of this subsampling estimator.

The next theorem verifies that, for independent data, the general results for convolved subsampling with previous block-based subsamples (Section~\ref{sec:mix}) also
hold when the convolution is based on the independent data subsampling estimator $\Lni$ using all   subsets of size $b$.

\begin{theoremnew1}\label{theorem7} Let $\{X_t\}$ be a sequence of independent (possibly non-iid) random variables.
 Given $\Lni$, let $\hatsigmaid$ and $\Cniid$ denote the corresponding   subsampling variance estimator  and  convolved subsampling estimator.
   Then,  Theorem~\ref{theorem4} holds under the notational convention that $\Ln\equiv \Lni$, $\hatsigma\equiv \hatsigmaid$, and $\Cn\equiv \Cniid$ and that subsample quantities in (\ref{eqn:Dn})-(\ref{eqn:cond})
are defined as $T_{b,i} = \tau_b[ t_b(Y_{b,i}) - t(P)]$, $i=1,\ldots,N_n\equiv {n \choose b}$.

Additionally, if the variables $\{X_t\}$ are iid, then Corollary~\ref{cor2} likewise holds.      \end{theoremnew1}

We may also draw some connections between convolved subsampling and the bootstrap
for sample means with independent data. Suppose independent variables $X_1,\ldots,X_n$ have common mean $\mu$
(e.g., as in Tukey's symmetric contamination model where observations may have different variances), from which we define $T_n \equiv \sqrt{n}(\bar{X}_n-\mu)$.
The convolved  estimator $\Cniid$ here has close parallels to the classic independent bootstrap of \citet{Efr}.
Namely, $\Cniid$ is the resampling distribution of $T_n^*\equiv \sqrt{n_1}(\bar{X}_{n_1}^* - \bar{X}_n)$ for a sample mean $\bar{X}_{n_1}^*$   of size $n_1=k_n b$ formed by averaging
 $k_n$ independent subsamples of size $b$, with each size $b$ subsample drawn uniformly and without replacement  from $\{X_i\}_{i=1}^n$; if the  subsamples of size $b$
are instead drawn with replacement from $\{X_i\}_{i=1}^n$, then  $T_n^*$ alternatively produces the independent bootstrap distribution, say  $\Cniidb$,  with a resample size $n_1= k_n b$.
Consequently,  the independent data version of  convolved subsampling $\Cniid$ does not exactly match the independent bootstrap. However, the following result for independent data shows that the  subsampling estimator $\Lni$ and its convolution $\Cniid$
 are valid in a broad non-iid context for  sample means, and the  differences between $\Cniid$
  and $\Cniidb$ are asymptotically negligible.

\begin{theoremnew1}\label{theorem8} Let $X_1,X_2,\ldots,$ denote a sequence of independent (possibly non-iid)   variables,
 with finite variances and common mean $\E X_t=\mu\in\mathbb{R}$. Define $X_{i,\mu}\equiv X_i-\mu$, $i\geq 1$.
% (i) Then, $T_n=\sqrt{n}(\bar{X}_n-\mu) \stackrel{d}{\rightarrow} N(0,\sigma^2)$ as $n\to \infty$  if,
% for each $\epsilon>0$ and some $\sigma^2>0$,
%\[
% \lim_{n \to \infty}  \frac{1}{n}\sum_{i=1}^n \E X_{i,c}^2I [|X_{i,c}|> \epsilon \sqrt{n}] =0\quad \mbox{and}\quad\lim_{n\to \infty} \left| \frac{1}{n}\sum_{j=1}^n \E X_{j,c}^2  -\sigma^2\right|=0.
%\]
%(ii)
As $n \to \infty$,  suppose
$b^{-1}+ b/n\to 0$  and that
\[
   \frac{1}{n}\sum_{i=1}^n \E X_{i,\mu}^2I [|X_{i,\mu}|> \epsilon \sqrt{b}] \to0\;\,\mbox{and}\;\,\max_{1 \leq i_1< i_2<\ldots < i_b \leq n} \left| \frac{1}{b}\sum_{j=1}^b \E X_{i_j,\mu}^2  -\sigma^2\right|\to 0
\]
  for each $\epsilon>0$ and some $\sigma^2>0$.
Then, as $n\to \infty$, $T_n=\sqrt{n}(\bar{X}_n-\mu) \stackrel{d}{\rightarrow} N(0,\sigma^2)$ along with $
\sup_{x \in \mathbb{R}} \lvert \Lni(x) - \Phi(x/ \sigma) \rvert \stackrel{p}{\rightarrow} 0$ and  $\hatsigmaid \stackrel{p}{\rightarrow} \sigma^2$.
Furthermore,  for any positive integer sequence $k_n$,
  \[
\sup_{x \in \mathbb{R}} |\Cniid(x) - \Phi(x/ \sigma)| \stackrel{p}{\rightarrow} 0 \ \ \ \ \text{and} \ \ \ \  d_2[\Cniid,\Cniidb]\stackrel{p}{\rightarrow} 0,
\]
where $d_2(\cdot,\cdot)$ denotes Mallow's metric between  distributions $\Cniid$ and  $\Cniidb$.
\end{theoremnew1}
 \noindent {\sc Remark 5}: Above $d_2(\cdot,\cdot)$ metricizes weak convergence (cf.~\citet{BiFr}), where, for distributions $F$ and $G$ on $\mathbb{R}$, $[d_2(F,G)]^2\equiv \inf\{ \E|X-Y|^2: X\sim F, Y\sim G \}$ with the infimum
over all random  pairs $(X,Y)$  with coordinate marginal distributions $F$ and $G$.\\

A  Lindeberg condition, defined by replacing $b$ with $n$ in   Theorem~\ref{theorem8} assumptions about second moment  limits,
suffices for $T_n=\sqrt{n}(\bar{X}_n-\mu) \stackrel{d}{\rightarrow} N(0,\sigma^2)$.    Hence,  Theorem~\ref{theorem8}  validates  subsampling, convolved subsampling and the bootstrap for the sample mean
 under a slightly stronger   condition than required for the CLT with independent data.   This  finding  also involves a
 weaker moment condition  than a classical   bootstrap result of \citet{Liu1988} for sample means of
 non-iid data (i.e., $\sup_{i \geq 1}\E |X_i|^{2+\delta}<\infty$ with $\delta>0$).   We add that, for
 approximating distributions of general statistics, \citet{PoRoWo}~(ch.~2.3) relate the subsampling
 estimator $\Lni$ itself to an independent-data version of the bootstrap with a resample size of $b<n$. For iid data,
 their Corollary~2.3.1  establishes the consistency of
the bootstrap estimator  from  $\Lni$ provided that $b^{-1}+b^2/n \to 0$  (i.e.,  to mitigate differences in resampling with or without replacement).
However, when specialized to sample means, their intended bootstrap estimator becomes
$\Cniidb$ with   convolution $k_n=1$, and the convolved subsampling estimator $\Cniid$ with $k_n=1$   matches $\Lni$ in the mean case.
For the sample mean, Theorem~\ref{theorem8}  consequently re-affirms the   equivalence
between   subsampling and bootstrap estimation
   with resample size $b<n$, but under a weaker  requirement $b/n\rightarrow 0$ and for non-iid data.

\section{Concluding remarks}
\label{sec:concl}
For approximating sampling distributions with normal limits, we have developed a general theory for estimators formed by the $k$-fold
self-convolution of subsampling  estimators.  Applied to time series,  convolved subsampling estimators have a close correspondence to the
block bootstrap, matching the latter for  sample means.  With more general statistics, convolved subsampling is
 not necessarily equivalent to the bootstrap, but instead provides a hybrid-type of resampling  for time series that has received recent
consideration in the literature (though without recognition as convolved subsampling).
As an advantage, convolved subsampling can be established under particularly weak process assumptions.
This facilitates, for example, proving and extending the consistency of the block bootstrap
under weaker assumptions than often previously considered.  The general theory here allows  convolved subsampling
estimation to be established from the consistency of the original subsampling estimator and its variance (i.e., a subsampling variance estimator),
where   the original subsampling estimator is often quite tractable  to verify and may be further applied to
 validating the subsampling variance estimator from general tools provided here.
It turns out that a consistent subsampling variance estimator is also necessary for consistent convolved subsampling, particularly with
increasing levels of convolution  which occurs in block bootstrap settings. Results here also broadly establish
convolved subsampling with general statistics for mixing and possibly non-stationary
processes, and we have examined the block bootstrap/convolved subsampling across several inference contexts
(mixing or linear or long-memory or spatial or independent processes), showing that the convolved subsampling
approach offers an alternative perspective for establishing
the bootstrap in the  sample mean case  under weak assumptions.
Based on emerging research directions with resampling  (cf.~Section~\ref{sec:convmore}),
 convolved subsampling and  results here may   provide tools  for  advancing future developments with dependent data.

\bibliographystyle{imsart-nameyear}
\bibliography{Lit}

\begin{appendix}

\section{Proofs of the main results}
\label{sec:proofs}
\newtheorem{lemma1}{Lemma}
\begin{lemma1} \label{lem1} Let $\{X_n\}$ be a sequence of real-valued random variables  such that $X_n \xrightarrow{d} X$ and $\var(X_n) \rightarrow \var(X)<\infty$ as $n \to \infty$. Then, $\E X_n \rightarrow \E X \in \mathbb{R}$ and $\E X_n^2 \rightarrow \E X^2<\infty$.
\end{lemma1}
\noindent {\it Proof.} Write $m_n =\E X_n$ and $Y_n=X_n - m_n$. As $\E Y_n^2$ is bounded, $\{Y_n\}$ is uniformly integrable and hence tight.  Because $\{X_n\}$ is also tight by  $X_n \xrightarrow{d} X$,
it holds that $m_n = X_n-Y_n$ is tight and so must be a bounded sequence.  Since
$\{Y_n\}$  and $\{m_n\}$ are uniformly integrable, the sum $\{X_n = Y_n +m_n\}$ is also uniformly
integrable and $\E X_n \rightarrow \E X$ follows. $\Box$  \\

\noindent \textbf{Proof of Theorem~\ref{theorem1}}. We show distributions converge in probability through Mallow's
metric $d_2(\cdot,\cdot)$; see Remark~5 and \citep{BiFr} (sec.~8).
For random variables $X,Y$ with $X\sim F$ and $Y\sim G$, we also denote $d_2(X,Y) \equiv d_2(F,G)$.

For $n\geq 1$, recall the random variable $Z_n^* \equiv k_n^{-1/2}\sum_{i=1}^{k_n}(Y_{n,i}^*-\mn)$ from (\ref{eqn:zn}) has  the
convolved distribution $\Cn$  based on iid variables $\{Y_{n,i}^*\}_{i=1}^{k_n}$ with distribution $\Ln$
and mean $\mn\equiv \int x d\Ln(x)$.  Let $Z_1,\ldots,Z_{k_n}$ be iid $N(0,\sigma^2)$ variables.
  Then,
\begin{eqnarray*}
[d_{2}(\Cn, \Phi(\cdot/\sigma))]^2 & =&   \left[d_{2}\left( Z_n^*, \frac{1}{\sqrt{k_n}}\sum_{i=1}^{k_n} Z_{i} \right)\right]^2 \\
& \leq &  \sum_{i=1}^{k_n} \frac{1}{k_n}\left[d_{2}\left(Y_{n,i}^*-\mn, Z_{i} \right)\right]^2 \\
 & = &     \left[d_{2}\left(Y_{n,1}^*, Z_{1} \right)\right]^2 - [\mn]^2
\end{eqnarray*}
holds by Lemmas 8.7-8.8 of \citep{BiFr} (for the inequality and the last equality, respectively).
By (\ref{eqn:1})-(\ref{eqn:2}) and $\hatsigma\stackrel{p}{\rightarrow }\sigma^2$, for any arbitrary subsequence
$\{n_j\} \subset \{n\}$, one may extract a further subsequence $\{\ell\equiv n_k\} \subset{n_j}$ such that
$Y_{\ell,1}^* \stackrel{d}{\rightarrow} Z_1$ and $\hatsigmal \rightarrow \sigma^2$ as $\ell\rightarrow \infty$ (a.s.$(P)$). By Lemma 1, this implies $\ml \rightarrow 0$ and $\int x^2 d \Ll(x)\rightarrow \sigma^2$ as $\ell\rightarrow \infty$ (a.s.$(P)$).  Together, $Y_{\ell,1}^* \stackrel{d}{\rightarrow} Z_1$ and $\int x^2 d \Ll(x)\rightarrow \sigma^2$ (where $Y_{\ell,1}^*\sim \Ll$) are equivalent to $[d_2(Y_{\ell,1}^*,Z_1)]^2 \rightarrow 0$ as $\ell\rightarrow \infty$ (cf.~Lemma 8.3, \citep{BiFr}).  Thus, $[d_{2}(\Cl, \Phi(\cdot/\sigma))]^2\rightarrow 0$, implying $\sup_{x\in\mathbb{R}}|\Cl(x)-\Phi(x/\sigma) |\rightarrow 0$ as $\ell \to \infty$ (a.s.$(P)$). As the subsequence $\{n_j\}$ was arbitrary, Theorem~\ref{theorem1} follows. $\Box$\\

\noindent \textbf{Proof of Proposition~\ref{prop1}}.  Recall that $\Cn(x)\equiv P_*(Z_n^* \leq x)$, $x\in\mathbb{R}$, is defined by the resampling distribution of $Z_n^*$ from (\ref{eqn:zn}), and define
 $\Cnu(x) \equiv P_*(Z_n^* + \mn \sqrt{k_n} \leq x)$, $x\in\mathbb{R}$,
where  $Z_n^* + \mn \sqrt{k_n} = \sum_{i=1}^{k_n} Y_{n,i^*}/\sqrt{k_n}$ with iid variables $\{Y_{n,i}^*\}_{i=1}^{k_n}$ drawn from
$\Ln$.
To show Proposition~\ref{prop1}, note that distributions $\Cn$  and $\Cnu$ can only match
  asymptotically   if and only if $\mn \sqrt{k_n} \stackrel{p}{\rightarrow} 0$, which is equivalent to  $\mn \equiv  \int x d \Ln(x) \stackrel{p}{\rightarrow} 0$  as $\{k_n\}$ is a bounded positive integer sequence.  Consequently,
 Proposition~\ref{prop1} follows by showing that, for bounded $\{k_n\}$,  (\ref{eqn:2}) holds   if and only if $\sup_{x\in \mathbb{R}}|\Cnu(x) - \Phi(x/\sigma)| \stackrel{p}{\rightarrow} 0$.

 Let $\phi_{n}(t)$ and $\tilde{\phi}_{n,k_n}(t) = [\phi_{n}(t/\sqrt{k_n})]^{k_n}$, $t\in \mathbb{R}$, denote the characteristic functions of $\Ln$ and $\Cnu$.  Suppose (\ref{eqn:2}) holds so that, for any subsequence $\{n_j\}\subset \{n\}$,
extract a further subsequence $\{\ell \equiv n_k\}\subset \{n_j\}$ such that $\sup_{x\in \mathbb{R}}|\Ll(x) - \Phi(x/\sigma)| \rightarrow 0$ as $\ell\rightarrow \infty$ a.s.$(P)$.  Then, for any given $T>0$,
\begin{align*}
\Delta_{\ell}(T)\equiv\sup_{|t|\leq T}|\phi_{\ell}(t) - e^{-t^2 \sigma^2/2}| \rightarrow 0
\end{align*}
 as $\ell\rightarrow \infty$ (a.s.$(P)$) by the Levy continuity theorem (cf.~\citep{KrishnaB.Athreya2006}, Theorem~10.3.1).
Fix $t_0 \in \mathbb{R}$ and set $T_0 = |t_0|$.  Then, using that
$|\prod_{i=1}^n w_i - \prod_{i=1}^n z_i| \leq \sum_{i=1}^n |w_i-z_i|$ for complex numbers $\{w_i,z_i\}_{i=1}^n$ with $|w_i|, |z_i|\leq 1$, we have
\[
  |\tilde{\phi}_{\ell,k_\ell}(t_0) - e^{-t_0^2 \sigma^2/2} | \leq k_\ell |\phi_{\ell}(t_0/\sqrt{k_n}) - e^{-t_0^2 \sigma^2/(2 k_\ell)}  |  \leq  k_\ell \Delta_{\ell}(T_0) \rightarrow 0
\]
as $\ell \rightarrow \infty$ (a.s.$(P)$) by $\sup_n k_n <\infty$.  Hence, for all $t_0\in\mathbb{R}$,
$\tilde{\phi}_{\ell,k_\ell}(t_0) \rightarrow e^{-t_0^2 \sigma^2/2}$ as $\ell \rightarrow \infty$, implying
$\sup_{x\in \mathbb{R}}|\Clu(x) - \Phi(x/\sigma)| \rightarrow 0$ (a.s.$(P)$).
As the subsequence $\{n_j\}$ was arbitrary, we have $\sup_{x\in \mathbb{R}}|\Cnu(x) - \Phi(x/\sigma)| \stackrel{p}{\rightarrow} 0$ by the equivalence of convergence in probability
to almost sure convergence along subsequences.

 Next suppose $\sup_{x\in \mathbb{R}}|\Cnu(x) - \Phi(x/\sigma)| \stackrel{p}{\rightarrow} 0$ so that,
  for any subsequence $\{n_j\}\subset \{n\}$,
extract a further subsequence $\{\ell\equiv n_k\}\subset \{n_j\}$ such
 that $\sup_{x\in \mathbb{R}}|\Clu(x) - \Phi(x/\sigma)| \rightarrow 0$ as $\ell\rightarrow \infty$ (a.s.$(P)$).
Fix $t_0$ and define $T = |t_0| \sup_n k_n$.  By  the Levy continuity theorem,
\begin{align*}
\sup_{|t|\leq T}\left| \tilde{\phi}_{\ell,k_\ell}(t) - e^{-t^2 \sigma^2/2}\right| =
\sup_{|t|\leq T}\left|[ {\phi}_{\ell}(t/\sqrt{k_\ell})]^{k_\ell} - [e^{-t^2 \sigma^2/(2 k_{\ell})}]^{k_{\ell}}\right| \rightarrow 0
\end{align*}
 as $\ell\rightarrow \infty$, implying $[\phi_\ell(t_0)/e^{-t_0^2 \sigma /2}]^{k_\ell}\rightarrow 1$
  (a.s.$(P)$) from bounded $\{k_\ell\}$.  As $1 \leq k_\ell \leq \sup_n k_n<\infty$,
  we then have $ \phi_\ell(t_0)\rightarrow e^{-t_0^2 \sigma /2} $
as $\ell \to \infty$ for any $t_0\in\mathbb{R}$, so that  $\sup_{x\in \mathbb{R}}|\Ll(x) - \Phi(x/\sigma)| \rightarrow  0$   (a.s.$(P)$).
As the subsequence $\{n_j\}$ was arbitrary,   $\sup_{x\in \mathbb{R}}|\Ln(x) - \Phi(x/\sigma)| \stackrel{p}{\rightarrow} 0$ or (\ref{eqn:2}) then follows, finishing the proof of Proposition~\ref{prop1}. $\Box$  \\

\noindent \textbf{Proof of Theorem~\ref{theorem2}}.
  For $n \geq 1$, again let $\{Y_{n,i}^*\}_{i=1}^{k_n}$ be iid with distribution $\Ln$ and mean $\mn \equiv \int x d\Ln(x)$, so that $Z_n^* \equiv k_n^{-1/2}\sum_{i=1}^{k_n}(Y_{n,i}^*-\mn)$ from (\ref{eqn:zn}) follows  the
convolved distribution $\Cn$.  For $\epsilon>0$, define quantities $\Delta_{1n}(\epsilon) \equiv   k_n P_* ( Y_{n,1}^* \geq  \sqrt{k_n} \epsilon) = k_n \int_{|x|\geq \sqrt{k_n} \epsilon} 1 d \Ln(x)$,
\[
\Delta_{2n}(\epsilon) \equiv     \sqrt{k_n} \int_{ |x|\geq \sqrt{k_n} \epsilon} |x| d\Ln(x), \qquad
\Delta_{3n}(\epsilon) \equiv  \int_{|x| \geq  \sqrt{k_n} \epsilon} |x|^2 d \Ln(x),
\]
 noting
\begin{align*}
\Delta_{1n}(\epsilon) \leq \epsilon\Delta_{2n}(\epsilon) \leq \epsilon^2\Delta_{3n}(\epsilon)
\end{align*}
 and, by assumption,  $\Delta_{3n}(\epsilon) \stackrel{p}{\rightarrow} 0$ for any $\epsilon>0$.
 For any $\{n_j\} \subset \{n\}$, extract a further subsequence $\{\ell\equiv n_k\}\subset\{n_j\}$ such that $\Delta_{3\ell}(1/m) \rightarrow 0 $  for any integer $m \geq 1$ as $\ell\rightarrow \infty$ (a.s.$(P)$),  implying that $\lim_{\ell\to \infty}\Delta_{j\ell}(\epsilon) =0 $ for any $\epsilon>0$ and $j=1,2,3$ (a.s.$(P)$). In
particular, as $|\tilde{\Delta}_{2\ell}(1)|\leq \Delta_{2\ell}(1)\to 0$ as $\ell\to\infty$ for $\tilde{\Delta}_{2\ell}(1) \equiv  \sqrt{k_\ell} \int_{ |x|\geq \sqrt{k_\ell} \epsilon} x d\Ll(x) $, note that $Z_\ell^*$ can have a normal $N(0,\sigma^2)$ limit law if and only if $Z_\ell^*-\tilde{\Delta}_{2\ell}(1)$ does as $\ell\rightarrow \infty$ (a.s.$(P)$).
  Also, from $\lim_{\ell \to \infty}\Delta_{1\ell}(\epsilon)=0$ for any $\epsilon > 0$, the array
$\{Y_{\ell,i}^*/\sqrt{k_\ell}\}_{i=1}^{k_\ell}$
is infinitesimal. Hence, by classical  CLT results   with independent, infinitesimal random variables (cf.~\citep{ChTe}, Theorem~3(ii), ch.~12.2),
\begin{align*}
 Z_\ell^*-\tilde{\Delta}_{2\ell}(1)=\frac{1}{\sqrt{k_{\ell}}}\sum_{i=1}^{k_{\ell}} \left [Y_{\ell,i}^*  - \int_{|x|<\sqrt{k_{\ell}}} x d \Ll(x) \right ]
\end{align*}
 will have a normal $N(0,\sigma^2)$ limit if and only if
\begin{align*}
  \Gamma_{\ell}(\epsilon) \equiv \int_{|x| < \sqrt{k_{\ell}}\epsilon} x^2 d \Ll(x) - \left( \int_{|x|<\sqrt{k_\ell} \epsilon} x d \Ll(x) \right)^2 \rightarrow \sigma^2
\end{align*}
holds for any  $\epsilon>0$ as $\ell\rightarrow \infty$ (a.s.$(P)$). But, for any $\epsilon>0$,
\begin{align*}
 & \ \lvert \Gamma_{\ell}(\epsilon) -\hatsigmal \rvert  \\
= & \ \left \lvert \int_{\lvert x  \rvert \geq \sqrt{k_{\ell}} \epsilon} x^2 d \Ll(x) +\left ( \int x d \Ll(x) \right )^2 - \left ( \int_{\lvert x \rvert < \sqrt{k_\ell} \epsilon} x d \Ll(x) \right )^2 \right \rvert \\
\leq & \ \Delta_{3\ell}(\epsilon) + 2 \Delta_{2\ell}(\epsilon) \int \lvert x \rvert d \Ll(x) \rightarrow 0
\end{align*}
by the almost surely boundedness of $\int \lvert x \rvert d \Ll(x)$ (as a consequence of $\Delta_{2\ell}(\epsilon) \rightarrow 0$). Thus, $\Gamma_{\ell}(\epsilon) \rightarrow \sigma^2$ for any $\epsilon>0$ if and only if $\hatsigmal \rightarrow \sigma^2$ as $\ell\rightarrow \infty$ (a.s.$(P)$).
As the subsequence
$\{n_j\}$ was arbitrary, the result now follows.$\Box$\\

\noindent \textbf{Proof of Corollary~\ref{cor1}}.   Suppose first $\hatsigma \stackrel{p}{\to}\sigma^2>0$ and $k_n\to\infty$.  Fix $\epsilon>0$.
Then, the normal convergence for $\Cn$ in  Corollary~\ref{cor1} will follow from Theorem~\ref{theorem2}
 by showing $\Upsilon_n(\epsilon) \equiv \int_{|x|\geq \epsilon\sqrt{k_n} } x^2 d\Ln(x) \stackrel{p}{\to}0$ under  any one of Conditions (C.1)-(C.4).    If Condition (C.1) holds, then for any subsequence $\{n_j\}\subset \{n\}$, extract a further subsequence $\{\ell\equiv n_k\}\subset\{n_j\}$ such that $\hatsigmal \rightarrow \sigma^2$ and $Y_\ell \stackrel{d}{\to} Y_0 $ as $\ell \to \infty$ (a.s.$(P)$), where $Y_\ell$ denotes a random variable with
  distribution $\Ll$ and $Y_0$ is a random variable with distribution $J_0$.  As $\hatsigmal$ and $\sigma^2$ are the variances of $Y_\ell$ and $Y_0$,   Lemma~\ref{lem1} yields $\int x^2 d \Ll(x) \to \E Y_0^2<\infty$ from which $\Upsilon_\ell(\epsilon) \to 0$ holds by $k_\ell \to \infty$ as $\ell \to \infty$ (a.s.$(P)$) (i.e., $\{Y_\ell^2: \ell \geq 1\}$ is uniformly integrable from $\int x^2 d \Ll(x) \to \E Y_0^2<\infty$ and  $Y_\ell \stackrel{d}{\to} Y_0 $).  As the subsequence $\{n_j\}$ was arbitrary,  $\Upsilon_n(\epsilon) \stackrel{p}{\to}0$  follows.
  If Condition (C.2) holds, then, for any $C>0$, Markov's inequality gives $P(\Upsilon_n(\epsilon) >  C) \leq P(\int  x^{2+\epsilon_0} d\Ln(x) >  (\epsilon\sqrt{k_n})^{\epsilon_0} C) \rightarrow 0$ as $n\to \infty$
  because $\int x^{2+\epsilon_0} d\Ln(x) = N_n^{-1} \sum_{i=1}^{N_n} [\tau_b(t_{n,b,i} -t_n)]^{2+\epsilon_0}=O_p(1)$ and $k_n\to \infty$.   If Condition (C.3)
   holds, then we use the inequality
   \[\Upsilon_n(\epsilon) \leq     2 \left(\frac{\tau_b T_n}{\tau_n}\right)^2 I(|\tau_b T_n/\tau_n|\geq 2^{-1}\epsilon\sqrt{k_n}) +\frac{2}{N_{n}}\sum_{i=1}T_{b,i}^2 I(|T_{b,i}|\geq 2^{-1}\epsilon\sqrt{k_n} ) \] to bound, for any $C>0$,
  \begin{eqnarray*}
&&  P(\Upsilon_n(\epsilon) >  C) \\& \leq &  P(|\tau_b T_n/\tau_n|\geq 2^{-1}\epsilon\sqrt{k_n})+\frac{4}{C} \frac{1}{N_n}\sum_{i=1}^{N_n} \E T_{b,i}^2 I(|T_{b,i}|\geq 2^{-1}\epsilon\sqrt{k_n} ) \to0
  \end{eqnarray*}
as $n\to \infty$ by $\sup_{b\geq 1}\sup_{1\leq i \leq N_n} \E T_{b,i}^2 I(|T_{b,i}|\geq 2^{-1}\epsilon\sqrt{k_n} ) \to 0$
and $|\tau_b T_n/\tau_n| = O_p(1)$.  Under Condition (C.4), the  above becomes
\[P(\Upsilon_n(\epsilon) >  C) \leq  C^{-1}4   \E T_{b}^2 I(|T_{b}|\geq 2^{-1}\epsilon\sqrt{k_n} ) + P(|\tau_b T_n/\tau_n|\geq 2^{-1}\epsilon\sqrt{k_n})\to 0\] with $|\tau_b T_n/\tau_n| = O_p(1)$ by $\sup_{n \geq 1}\E T_n^2 <\infty$ and $\tau_b/\tau_n=O(1)$.

Next consider the converse of Corollary~\ref{cor1}.   For  $A_n \equiv \max_{1\leq i \leq N_n} |\tau_b[t_{n,b,i}-t_n]|$, $n\geq 1$,  pick an increasing integer sequence $k_{n+1}>k_{n} \geq 1$ such that
$P(A_n \geq \sqrt{k_n} ) < 2^{-n}$ for each $n\geq 1$, where $\{\tau_b[t_{n,b,i}-t_n]\}_{i=1}^{N_n}$
are the subsample statistics defining $\Ln$ in (\ref{eqn:sub}).
 By the Borel-Cantelli lemma, $A_n<\sqrt{k_n}$ holds eventually for large $n$ (a.s.$(P)$).  Next recall that $\Cn$ corresponds to the distribution of $Z_n^* \equiv k_n^{-1/2}\sum_{i=1}^{k_n} (Y_{n,i}^*-\mn)$ from (\ref{eqn:zn}), where $\{Y_{n,i}^*\}_{i=1}^{k_n}$ are iid variables following $\Ln$.
 Fix an arbitrary integer $m \geq 1$.  As $k_n\to \infty$,
it holds that $ P^*( |Y_{n,1}^*| > m^{-1}\sqrt{k_n}) \leq \Ln(-m^{-1}\sqrt{k_n})+ 1 - \Ln(m^{-1}\sqrt{k_n}) \stackrel{p}{\to}0$ under any one of the Conditions  (C.1)-(C.4), which can be established similarly to arguments above.
 As $\Cn$ converges by assumption for $\{k_n\}$,
then for any subsequence $\{n_j\}\subset \{n\}$, we may extract a further subsequence $\{\ell\equiv n_k\}\subset\{n_j\}$ such that
$\sup_{x\in\mathbb{R}}|\Cl(x)-\Phi(x/\sigma)|\to0$  holds, along with $ P^*( |Y_{\ell,1}^*| > m^{-1}\sqrt{k_\ell})\to 0$ for any $m \geq 1$, as $\ell \to \infty$ (a.s.$(P)$).  As the row-wise independent array $\{Y_{\ell,i}^*/\sqrt{k_\ell}\}_{i=1}^{k_\ell}$, $\ell \geq 1$, is infinitesimal and
$Z_\ell^*$ has a normal $N(0,\sigma^2)$ limit, it follows that
\[
\lim_{\ell \to \infty} \left[ \int_{|x| < \sqrt{k_\ell}} x^2 d\Ll(x) -\left( \int_{|x| < \sqrt{k_\ell}} x d\Ll(x)\right)^2 \right]=\sigma^2
\]
(a.s.$(P)$) by classical convergence results to normal laws (\citep{ChTe}, ch.~12.2, Theorems~2-3).  However,
 $\int_{|x| < \sqrt{k_\ell}} x^2 d\Ll(x) = \int  x^2 d\Ll(x)$ and  $\int_{|x| < \sqrt{k_\ell}} x  d\Ll(x) =\int x  d\Ll(x)=\ml$
 eventually for large $\ell$ (a.s.$(P)$) as $A_\ell < \sqrt{k_\ell}$ eventually.  Hence, \[\lim_{\ell \to \infty}  \hatsigmal =
   \lim_{\ell \to \infty} [\int x^2 d\Ll(x) - ( \ml )^2 ]=\sigma^2.\]  Now $\hatsigma \stackrel{p}{\rightarrow}\sigma^2$
   follows in Corollary~\ref{cor1} as $\{n_j\}$ was arbitrary.
$\Box$\\

\noindent\textbf{Proof of Theorem~\ref{theorem3}}.  We first establish $\hatsigma \stackrel{p}{\rightarrow}\sigma^2$  assuming (\ref{eqn:cond}). Let $\epsilon > 0$ and $\delta>0$, set
$Z \sim N(0,\sigma^2)$, and define $\Lt(x) \equiv N_n^{-1} \sum_{i=1}^{N_n} I[T_{b,i}\leq x] = \Ln(x + \tau_b(t_n-t(P)))$, $x\in\mathbb{R}$,
as the empirical distribution of the subsample copies $T_{b,i} = \tau_b[t_b(X_i,\ldots,X_{i+b-1}) -t(P)]$,
$i = 1,\ldots,N_n\equiv n-b+1$ of length $b \equiv b_n \in [1, n)$. By (\ref{eqn:cond}) and noting
$\int_{|x|>m} x^2 d\Lt(x) = N_n^{-1} \sum_{i=1}^{N_n} T_{b,i}^2I[|T_{b,i}|>m]$ for $m>0$, choose and fix integer $m \geq 1$ such that $\E Z^2 I[|Z|>m]<\epsilon/3$ and $P(\int_{|x|>m} x^2 d\Lt(x) >\epsilon/3)<\delta$ for all $n \geq m$.  For any $\{n_j\}\subset\{n\}$, extract a further subsequence
$\{\ell\equiv n_k\}\subset\{n_j\}$ such that $\sup_{x\in\mathbb{R}}|\Ll(x)-\Phi(x/\sigma)|\rightarrow 0$
and $\tau_{b_\ell}[t_\ell-t(P)]\rightarrow 0$ as $\ell\rightarrow \infty$ (a.s.$(P)$) by (\ref{eqn:2}) and $T_n \equiv \tau_n[t_n-t(P)]=o_p(\tau_n/\tau_b)$.  Hence, as $\ell\rightarrow \infty$ (a.s.$(P)$), these two limits imply $\sup_{x\in\mathbb{R}}|\lt(x)-\Phi(x/\sigma)|\rightarrow 0$ holds as well as, by the Dominated Convergence Theorem (DCT), $\int_{|x|\leq m} x^2 \lt(x) \to \E Z^2 I[|Z|\leq m]$.  As the subsequence
$\{n_j\}$ was arbitrary,  $\sup_{x\in\mathbb{R}}|\Lt(x)-\Phi(x/\sigma)|\stackrel{p}{\rightarrow} 0$
and  $\int_{|x|\leq m} x^2 d\Lt(x) \stackrel{p}{\rightarrow} \E Z^2 I[|Z|\leq m]$ hold as $n\to \infty$.
Then, we may bound
\begin{eqnarray*}
&&\lim_{n\to\infty} P\left( \left|\int x^2 d\Lt(x) -\sigma^2\right|>\epsilon\right)\\
&\leq &\lim_{n\to\infty} P\left( \left|\int_{|x|>m} x^2 d\Lt(x)\right|>\epsilon/3\right) \\
&& \qquad+\lim_{n\to\infty} P\left( \left|\int_{|x|\leq m} x^2 d\Lt(x) - \E Z^2 I[|Z|\leq m]\right|>\epsilon/3\right)\\
&\leq& \delta
\end{eqnarray*}
so that, as $\epsilon,\delta>0$ were arbitrary, $\int x^2 d\Lt(x) \stackrel{p}{\rightarrow}
\sigma^2$ follows directly, implying also $\int x  d\Lt(x) \stackrel{p}{\rightarrow}
0$ from   $\sup_{x\in\mathbb{R}}|\Lt(x)-\Phi(x/\sigma)|\stackrel{p}{\rightarrow} 0$.
Now $\hatsigma \stackrel{p}{\rightarrow}
\sigma^2$ follows by the expansion
\begin{equation}
\label{eqn:hatsig}
\hatsigma = \int x^2 d\Lt(x)-\left( \int x  d\Lt(x)\right)^2.
\end{equation}

We next show (\ref{eqn:cond}) holds if $\hatsigma \stackrel{p}{\rightarrow}
\sigma^2$. If possible, suppose (\ref{eqn:cond}) does not hold so that,
for some $\epsilon > 0$ and some $\delta > 0$, there exists an integer subsequence pair $\{(m_j,n_j)\}_{j \geq 1}$, with
$\{m_j\}\subset \{m\}$ and $n_j \geq  m_j > n_{j-1} \geq m_{j-1}$ for each $j > 1$, such that
\begin{equation}\label{eqn:bad}
P\left(
\int_{|x|>m_j} x^2 d\tilde{S}_{n_j}(x) >\epsilon\right) >\delta\end{equation}
for all $j\geq 1$, using again that
$\int_{|x|>m} x^2 d\Lt(x) = N_n^{-1} \sum_{i=1}^{N_n} T_{b,i}^2I[|T_{b,i}|>m]$, $n,m\geq 1$.
However, from  $\hatsigma \stackrel{p}{\rightarrow}
\sigma^2$ and (\ref{eqn:2}) and $T_n = o_p(\tau_n/\tau_b)$, we may extract a further subsequence
$\{\ell \equiv n_k\}\subset\{n_j\}$ such that $\sup_{x\in\mathbb{R}}|\lt(x)-\Phi(x/\sigma)| \rightarrow  0$ (in the same manner
as above) and $\hatsigmal  \rightarrow
\sigma^2$ as $\ell \to \infty$ (a.s.($P$)). From this, the subsampling variance
expansion (\ref{eqn:hatsig}) and Lemma~1, the sequence of random variables $\{\tilde{Y}_\ell^2\}_{\ell \geq 1}$, where $\tilde{Y}_\ell$
has distribution $\lt$ and variance $\hatsigmal$, must be uniformly integrable (a.s.($P$)).
That
is, $\lim_{m\to \infty} \sup_{\ell \geq 1} \int_{|x| >m }x^2 d \lt(x) =0$ (a.s.($P$)),
 which creates a contradiction of (\ref{eqn:bad}): letting
$\tilde{m}_1, \tilde{m}_2, \ldots$ denote the subsequence $\{m_k\}$  of $\{m_j\}$  paired with $\{\ell \equiv n_k\}\subset\{n_j\}$,
\[
0 < \delta < \sup_{\ell \geq 1} P\left( \int_{|x|>\tilde{m}_j} x^2 d\lt(x) > \epsilon\right)\leq
P\left(  \sup_{\ell \geq 1}  \int_{|x|>\tilde{m}_j} x^2 d\lt(x) > \epsilon\right)\rightarrow 0
\]
as $j \to \infty$.

For part(ii) of Theorem~\ref{theorem3}, if $Y_b^2$, $b\equiv b_n \geq 1$, is uniformly integrable (for $Y_b$ with
distribution $D_{n,b}$ from (\ref{eqn:Dn})), then for any $\epsilon> 0$ and $m \geq  1$,
\begin{eqnarray*}
\sup_{n \geq m} P\left( \int_{|x|>m} x^2  d\Lt(x)  >\epsilon\right) &\leq& \frac{1}{\epsilon} \sup_{n \geq 1} \E \left( \int_{|x|>m} x^2  d\Lt(x)\right) \\&=& \frac{1}{\epsilon} \sup_{n \geq 1} \int_{|x|>m} x^2 d D_{n,b}(x)\rightarrow 0
\end{eqnarray*}
  as $m\to \infty$, showing (\ref{eqn:cond}). This also establishes (\ref{eqn:cond}) under the assumptions of part(ii), where
$P(T_b \leq x) = D_{n,b}(x)$, $x \in \mathbb{R}$ and $T_b^2$ is uniformly integrable.   In part(ii) of Theorem~\ref{theorem3}, note that (\ref{eqn:2}) and $T_n =o_p(\tau_n/\tau_b)$ entail that $\Lt(x)\stackrel{p}{\to} \Phi(x/\sigma)$
for each $x\in\mathbb{R}$.  Hence, by the DCT applied to $|\Lt(x)|\leq 1$, then $ D_{n,b}(x) = \E\Lt(x) \to  \Phi(x/\sigma)$
holds for any $x\in\mathbb{R}$ so that $Y_b \stackrel{d}{\to} N(0,\sigma^2)$
 follows. Consequently,
$ \E Y_b^2 \to \sigma^2$ is equivalent to $Y_b^2$, $b\equiv b_n \geq 1$, being uniformly integrable.
$\Box$\\

\noindent \textbf{Proof of Theorem~\ref{theorem4}}.
By the assumptions, (\ref{eqn:2}) follows (i.e., $\Ln$ converges to a normal
limit) by Theorem~4.2.1 of \citep{PoRoWo}. Then, assumption (\ref{eqn:cond}) with (\ref{eqn:2}) gives $\hatsigma \stackrel{p}{\rightarrow} \sigma^2$ by Theorem~\ref{theorem3} and the convergence of $\Cn$ follows from Theorem~\ref{theorem1}. $\Box$\\

\noindent \textbf{Proof of Theorem~\ref{theorem5}}.  We first establish a CLT for $T_n = \sqrt{n}(\bar{X}_n -\mu) =
n^{-1/2}\sum_{t=1}^n (X_t - \E X_t)$ using results from
  \citep{KrishnaB.Athreya2006}~(ch.~16). For integers $n \geq 1$, $i \geq 1$, and real $M>1$, define sum quantities $T_{n,i} \equiv n^{-1/2} \sum_{t=i}^{i+n-1} (X_t- \E X_t)$ as well as truncated versions
  \begin{eqnarray*}
  T_{n,i}^{(1)}(M) &\equiv &   n^{-1/2} \sum_{t=i}^{i+n-1} (X_t I(|X_t| \leq M)- \E X_t I(|X_t| \leq M))\\
  T_{n,i}^{(2)}(M) &\equiv &   n^{-1/2} \sum_{t=i}^{i+n-1} (X_t I(|X_t| > M)- \E X_t I(|X_t|> M)) = T_{n,i}-T_{n,i}^{(1)}(M)
  \end{eqnarray*}
  for any $i,n,M$. By assumption,
  \begin{equation}
  \label{eqn:var2}
  \sup_{i \geq 1} |\var(T_{n,i})-\sigma^2| \rightarrow 0 \quad \mbox{as $n\to\infty$}
  \end{equation}
  holds and  \citep{KrishnaB.Athreya2006}~(p. 526) show that, under the mixing and moment
assumptions, there exists some $C > 0$ (not depending on $M$ or $n$) such that
  \begin{equation}
  \label{eqn:var3}
  \sup_{i \geq 1} \E [T_{n,i}^{(2)}(M)]^2 \leq C \left( M^{-3 \delta /4} + \sum_{k=\lfloor M^{\delta/4}\rfloor}^\infty \alpha(k)^{\delta/(2+\delta)}\right) \equiv \Lambda(M)
  \end{equation}
 for all $M>1$ and $n \geq 1$; note $\lim_{M\to \infty} \Lambda(M)=0$ also holds by the mixing assumptions. The proof of \citep{KrishnaB.Athreya2006}~(Theorem~16.3.2) provides a CLT (\ref{eqn:1}) for $T_n \equiv T_{n,1}$
assuming bounded random variables $\{X_t\}$, but the same arguments hold immediately for
$T^{(1)}_{n,1}(\log n)$ (i.e., variables truncated at $\log n$) provided that $\lim_{n\to \infty} \sup_{i \geq 1} |\var(T^{(1)}_{n,i}(\log n))-\sigma^2|=0$.  The latter follows by
(\ref{eqn:var2})-(\ref{eqn:var3}) here, so that we have $T^{(1)}_{n,1}(\log n)\stackrel{d}{\rightarrow} N(0,\sigma^2)$ as $n\to \infty$.  Also, for $\imath=\sqrt{-1}$ and $t\in\mathbb{R}$, if $\phi_{b,i}^{(1)}(t) \equiv \E e^{\imath t T_{b,i}^{(1)}(\log b)}$ denotes the   characteristic function of $T_{b,i}^{(1)}(\log b)$ for $i=1,\ldots,N_n \equiv n-b+1$, the   same proof of Athreya and Lahiri (2006) (Theorem~16.3.2) shows
 \begin{equation}
  \label{eqn:var4}
  \max_{1 \leq i \leq N_n} \left| \phi_{b,i}^{(1)}(t) - e^{-t^2 \sigma^2/2}\right|\to 0,
  \end{equation}
for each $t\in\mathbb{R}$ as $n\to \infty$ using (\ref{eqn:var2})-(\ref{eqn:var3}) with $b^{-1}+b/n\rightarrow 0$.

Now from $T^{(1)}_{n,1}(\log n)\stackrel{d}{\rightarrow} N(0,\sigma^2)$ and $T^{(2)}_{n,1}(\log n)\stackrel{p}{\rightarrow} 0$  as $n\to \infty$, where the latter follows
from $\E[ T^{(2)}_{n,1}(\log n)]^2\leq \Lambda(\log n)\to 0$ under (\ref{eqn:var3}), we obtain
 $T_n= T^{(1)}_{n,1}(\log n)+T^{(2)}_{n,1}(\log n)\stackrel{d}{\rightarrow} N(0,\sigma^2)$ in Theorem~\ref{theorem5} by Slutsky's theorem.

 Furthermore,  if $\phi_{b,i}(t) \equiv \E e^{\imath t T_{b,i}}$ denotes the   characteristic function of $T_{b,i}= T_{b,i}^{(1)}(\log b) + T_{b,i}^{(2)}(\log b)$, $1 \leq i \leq N_n$, then
\begin{eqnarray*}
&&\max_{1 \leq i \leq N_n} \left| \phi_{b,i} (t) - e^{-t^2 \sigma^2/2}\right|\\
&\leq& \max_{1 \leq i \leq N_n} \left| \phi_{b,i}^{(1)}(t) - e^{-t^2 \sigma^2/2}\right|+
\max_{1 \leq i \leq N_n} \left| \phi_{b,i} (t) - \phi_{b,i}^{(1)}(t)\right|
\to 0
\end{eqnarray*}
for each $t \in \mathbb{R}$ as $n \to 1$ by (\ref{eqn:var4}) and
\[
\max_{1 \leq i \leq N_n} \left| \phi_{b,i} (t) - \phi_{b,i}^{(1)}(t)\right| \leq \max_{1 \leq i \leq N_n}
\E\left| e^{\imath t T_{b,i}^{(2)}(\log b)}-1\right|\leq |t| [\Lambda(\log b)]^{1/2}\to 0
\]
 from (\ref{eqn:var3}) along with $|e^{\imath (x+y)} - e^{\imath y}| =|e^{\imath x}-1|\leq |x|$ for $x,y\in\mathbb{R}$.  From this, we obtain that
if $Y_b$ denotes a random variable with distribution function $D_{n,b}$ from (\ref{eqn:Dn}), then the characteristic
function of $Y_b$  satisfies $\E e^{\imath t Y_b}\equiv N_{n}^{-1}\sum_{i=1}^{N_n} \phi_{b,i}(t) \to e^{-t^2 \sigma^2/2}$ for each
$t \in \mathbb{R}$ as $n\to \infty$.    Hence, $Y_b \stackrel{d}{\rightarrow}N(0,\sigma^2)$
as $n \to \infty$  which further implies $\{Y_b^2: b \geq 1\}$  is uniformly integrable
(cf.~Lemma 1) as the second moment of $Y_b$ here is $\E Y_b^2 \equiv N_n^{-1}\sum_{i=1}^{N_n} \var(T_{b,i}) \to \sigma^2$ by (\ref{eqn:var2}).   Now Theorem~\ref{theorem5}
follows from Theorem~\ref{theorem4} as $T_n = O_p(1) = o_p((n/b)^{1/2})$  and $Y_b \stackrel{d}{\rightarrow}N(0,\sigma^2)$  hold and, by
Theorem~\ref{theorem3}(ii), (\ref{eqn:cond}) does as well. $\Box$\\

\noindent \textbf{Proof of Corollary~\ref{cor3}}. Recall here $T_n \equiv \sqrt{n} (\bar{X}_n-t(P))$ for $t(P) =M(\mu)$  based on
$\mu(t)=   \E X_t$, $t \in \mathbb{Z}$, and the subsample estimator $\Ln$ is the   empirical distribution (\ref{eqn:sub}) of $\sqrt{b}(\bar{X}_{b,i}-\bar{X}_n)$ for $\bar{X}_{b,i} \equiv \sum_{t=i}^{i+b-1}X_t/b$, $i=1,\ldots,N_n\equiv n-b+1$.

Consider the time series $\{Y_t\}$, defined by $Y_t = X_t - \mu(t)$, having mean zero $\mu_Y= 0$.
From a sample $Y_1,\ldots,Y_n$, write $\TnY \equiv \sqrt{n}(\bar{Y}_n-\mu_Y)$ based on the sample mean $\bar{Y}_n$, and let $\hatsigmaY$ and $\LnY$ denote the subsample variance and distribution estimators for $\TnY$ as derived from the subsample  quantities $b^{1/2}(\bar{Y}_{b,i}-\bar{Y}_n)$  for $\bar{Y}_{b,i} \equiv \sum_{t=i}^{i+b-1}Y_t/b$, $1 \leq i \leq N_n$. As $\{Y_t\}$  is an APC strongly mixing time series, Lemma A.6 of \citep{Syn} yields
\[
\sup_{i \geq 1} \left| \var\left( n^{-1/2}\sum_{t=i}^{i+n-1}Y_t\right)-\sigma^2\right|\to 0
\]
as $n\to \infty$ for some $\sigma>0$.
The assumptions of Theorem~\ref{theorem5}  then hold for $\{Y_t\}$ so
that
\begin{equation}
\label{eqn:apc1}
\TnY \stackrel{d}{\rightarrow} N(0,\sigma^2),\qquad \sup_{x\in\mathbb{R}}|\LnY(x) - \Phi(x/\sigma)|\stackrel{p}{\rightarrow}0,\qquad \hatsigmaY \stackrel{p}{\rightarrow} \sigma^2,
\end{equation}
as $n\to \infty$ with $b^{-1}+b/n\to 0$. Using (\ref{InequalAlmostPeriodic}), it holds that
\begin{equation}
\label{eqn:apc2}
\sup_{i \geq 1} \left| \frac{1}{n} \sum_{t=i}^{i+n-1} [\mu(t)-M(\mu)]\right|\leq \frac{C}{n},\qquad
\sup_{i \geq 1} \left| \bar{X}_{b,i}-\bar{Y}_{b,i}-M(\mu) \right|\leq \frac{C}{b}
\end{equation}
for all $b,n\geq 1$ with some $C>0$ (not depending on $n,b$).  By (\ref{eqn:apc1})-(\ref{eqn:apc2}), the limit distribution of $T_n$ follows as $T_n \stackrel{d}{\rightarrow}N(0,\sigma^2)$ by $T_n - \TnY = n^{-1/2}\sum_{t=1}^n [\mu(t)-M(\mu)] = O(n^{-1/2})$.   Likewise, by (\ref{eqn:apc2}), the difference of subsample statistics
\[
d_n \equiv \max_{1 \leq i \leq N_n} |\sqrt{b}(\bar{X}_{b,i}-\bar{X}_n)-\sqrt{b}(\bar{Y}_{b,i}-\bar{Y}_n)|\leq C b^{1/2}(b^{-1}+n^{-1})\to 0
\]
as $n\to \infty$.  Consequently, using that $\LnY(x-d_n) \leq \Ln(x) \leq \LnY(x+d_n)$ holds for all $x\in\mathbb{R}$, we find that
\begin{eqnarray*}
&&\sup_{x\in\mathbb{R}} |\Ln(x)-\Phi(x/\sigma)| \\&\leq&  \sup_{x\in\mathbb{R}} |\LnY(x)-\Phi(x/\sigma)|  +
\sup_{x\in\mathbb{R}} |\Phi((x+d_n)/\sigma)-\Phi(x/\sigma)| \stackrel{p}{\rightarrow}0
\end{eqnarray*}
by (\ref{eqn:apc1}) and the continuity of $\Phi(\cdot)$.  Hence, (\ref{eqn:2}) holds or $\Ln$ is consistent. Finally, by (\ref{eqn:apc1}), Theorem~\ref{theorem3}(i) yields
\begin{equation}
\label{eqn:apc3}
\lim_{m\to \infty} \Delta_m(\epsilon)=0,\qquad \Delta_m(\epsilon)\equiv \sup_{n \geq m} P\left( N_n^{-1} \sum_{i=1}^{N_n} b [\bar{Y}_{b,i}]^2 I(\sqrt{b}|\bar{Y}_{b,i}| >m) > \epsilon\right)
\end{equation}
for each $\epsilon>0$.  Fixing $\epsilon>0$ and using (\ref{eqn:apc2}), we have
\begin{eqnarray*}
&&\sup_{n \geq m} P\left( N_n^{-1} \sum_{i=1}^{N_n} b [\bar{X}_{b,i}-M(\mu)]^2 I(\sqrt{b}|\bar{X}_{b,i}-M(\mu)| >m) > \epsilon\right)\\
&\leq & \sup_{n \geq m} P\left(b^{-1} C^2 +  N_n^{-1} \sum_{i=1}^{N_n} b [\bar{Y}_{b,i}]^2 I(\sqrt{b}|\bar{Y}_{b,i}| >m -b^{-1/2}C) > \epsilon/2\right)\\&\leq& \Delta_{m-1}(\epsilon/4)
\end{eqnarray*}
where the last inequality follows for any large $m$ such that $b^{-1/2}C<1$  and $b^{-1}C^2<\epsilon/4$ hold
for $n \geq m - 1$  based on $C > 0$ in (\ref{eqn:apc2}). By this, (\ref{eqn:apc3}) and (\ref{eqn:2}) (i.e., $\Ln$ consistency),
$\hatsigma \stackrel{p}{\rightarrow}\sigma^2 $  holds by Theorem~\ref{theorem3}(i) and then $\sup_{x\in\mathbb{R}} |\Cn(x) - \Phi(x/\sigma)| \stackrel{p}{\rightarrow}0$   follows by
Theorem~\ref{theorem1}. $\Box$\\

\noindent \textbf{Proof of Corollary~\ref{cor4}}. By (\ref{eqn:1}) and (\ref{eqn:var}), $\{T_n^2:n \geq 1\}$   is uniformly integrable where
$T_n \equiv n^{\alpha/2}(\bar{X}_n-\mu)$.    This implies (\ref{eqn:cond}) by Theorem~\ref{theorem3}(iii). Now Corollary~\ref{cor4} follows from
Corollary~1 under (\ref{eqn:2}). $\Box$\\

\noindent \textbf{Proof of Corollary~\ref{cor5}}. This follows as a special case of Corollary~\ref{cor4}. $\Box$\\

 \noindent\textbf{Proof of Theorem~\ref{theorem6}}.    Let $\Lnt$ denote the subsampling estimator found by using $T_{b,\bm{i}}\equiv \sqrt{n_b}(\bar{X}_{b,\bm{i}}-\mu)$, $\bm{i}\in \mathcal{I}_n$, in (\ref{eqn:sub}).  By the boundary conditions on $R_0$ and
$\mathcal{D}_0$, the number $N_n=|\mathcal{I}_n|$ of subsamples, the sample size $n$ and the subsample size $n_b$ satisfy $N_n/\mathrm{vol}(R_n)\rightarrow 1$, $n/\mathrm{vol}(R_n)\rightarrow 1$ and $n_b/[b^d \mathrm{vol}(\mathcal{D}_0)]\rightarrow 1$ where $\mathrm{vol}(R_n)=\lambda_n^d \mathrm{vol}(R_0)$; see \citep{Lah2}, ch.~12.
By \citep{Lah2}~(Proposition~4.1 and Theorem~4.3), $T_n \stackrel{d}{\rightarrow }N(0,\sigma^2)$
and $\var(T_n)\rightarrow \sigma^2$ as $n\to \infty$; the same is true replacing $T_n$ with  $T_{b,\bm{0}}$ for the origin $\bm{0}\in\mathbb{R}^d \cap \mathcal{I}_n$.
Hence,  $\{T_{b,\bm{i}}:\bm{i}\in \mathcal{I}_n\}$ are uniformly integrable by stationarity and, by Theorem~\ref{theorem3}(iii),
the consistency of $\hatsigma$ (and of $\Cn$ from Theorem~\ref{theorem1}) will follow by showing (\ref{eqn:2})
or normality of $\Ln$.  Since $\Lnt(x+ T_n[\sqrt{n_b/n}]) = \Ln(x)$ for each $x\in\mathbb{R}$
and $T_n[\sqrt{n_b/n}]=o_p(1)$ by $b/\lambda_n\rightarrow 0$, it suffices to show
$ |\Lnt(x)-\Phi(x/\sigma)|\stackrel{p}{\rightarrow} 0$ for any $x$.  To this end, as $\E \Lnt(x) = P(T_{b,\bm{0}} \leq x )\rightarrow \Phi(x/\sigma)$ for $x\in \mathbb{R}$, it remains to show $\var(\Lnt(x))=o(1)$.
Let $ \sigma_x(\bm{i}) \equiv \mathrm{Cov}[I(T_{b,\bm{0}} \leq x), I(T_{b,\bm{i}} \leq x)]$, $\bm{i}\in \mathbb{Z}^d$. For any $x\in \mathbb{R}$, note $|\sigma_x(\bm{i})|\leq 1$ holds for $\bm{i}\in\mathbb{Z}^d$, while $|\sigma_x(\bm{i})|\leq 4 \alpha( |\bm{i}|_\infty-b;b)$
for $ \bm{i}\in\mathbb{Z}^d$ with $|\bm{i}|_\infty>b$ (cf.~Corollary 16.2.4(ii), \citep{KrishnaB.Athreya2006}).
By stationarity, the covariance bounds, and $|\{\bm{i}\in\mathbb{Z}^d: |\bm{i}|_{\infty}=k\}|\leq 2 d (2k+1)^{d-1}$ for $k \geq 0$, we have (for a generic  $C>0$)
\begin{eqnarray*}
\var(\Lnt(x))
&\leq &\frac{C}{N_n} \sum_{k=0}^{2b} (k+1)^{d-1} +  \frac{C}{N_n} \sum_{k=2b+1}^{\lambda_n} (k+1)^{d-1}
\alpha(k-b; b)\\ &\leq& O\left( \frac{b^d}{\lambda_n^d}\right ) + \frac{C}{N_n} \sum_{k=b+1}^{\infty} k^{d-1-\tau_1+\tau_2}=o(1)
\end{eqnarray*}
as $\tau_1-d-\tau_2>0$, $\tau_2 \geq 0$, and $b/\lambda_n + N_n^{-1}\rightarrow 0$.  $\Box$\\

 \noindent\textbf{Proof of Theorem~\ref{theorem7}}.   The proof essentially matches the one for Theorem~\ref{theorem4} except that, in place of mixing, one uses independence to show $\var(\Lni(x)) \leq \lfloor n/b\rfloor^{-1} 4^{-1}\rightarrow 0$ holds for any $x\in\mathbb{R}$, as in Remark~2.2.3 of \citep{PoRoWo} (applying to independent data). $\Box$\\

 \noindent\textbf{Proof of Theorem~\ref{theorem8}}.  Conditions of    Theorem~\ref{theorem8}(ii) imply those of  Theorem~\ref{theorem8}(i) so that
 $\sqrt{n}(\bar{X}_n -\mu)\stackrel{d}{\to} N(0,\sigma^2)$ holds.  By this and $b^{-1}+b/n\rightarrow 0$, the convergence results
 for $\Lni$, $\hatsigmaid$ and $\Cniid$ in Theorem~\ref{theorem8}(ii) will follow from Theorem~\ref{theorem7} by showing that (\ref{eqn:cond}) holds
  and that $Y_b \stackrel{d}{\to} N(0,\sigma^2)$ for a random variable $Y_b$ with distribution $D_{n,b}$ from (\ref{eqn:Dn}). In the   case here, note
 (\ref{eqn:Dn})-(\ref{eqn:cond}) are defined by subsample quantities $T_{b,i}$, $i=1,\ldots,N_n\equiv {n \choose b}$, of form
 $(X_{i_1,\mu}+\cdots +X_{i_b,\mu})/\sqrt{b}$ for some $1 \leq i_1<i_2<\cdots <i_b \leq n$ with $X_{i,\mu}=X_i-\mu$, $i \geq 1$.   Furthermore, by Theorem~\ref{theorem3}(ii),  (\ref{eqn:cond})    follows by $Y_b \stackrel{d}{\to} N(0,\sigma^2)$ because
 \[
 \E Y_b^2 = \frac{1}{N_n} \sum_{i=1}^{N_n} \E T_{b,i}^2 = \frac{1}{{n\choose b}} \sum_{1 \leq i_1<i_2<\cdots <i_b \leq n} \left( \frac{1}{b} \sum_{j=1}^b \sigma_{i_j}^2 \right)= \frac{1}{n}\sum_{i=1}^n \sigma_i^2 \equiv v_n \to \sigma^2
 \]
 as $n\to \infty$ from the assumed conditions, where we denote $\sigma_i^2 \equiv \var(X_i)=\E X_{i,\mu}^2$, $i \geq 1$.  Note that above, and in the following, we use that each index $1 \leq j \leq n$ appears exactly ${n-1 \choose b-1}$ times in the sum over subsample sets $1 \leq i_1<i_2<\cdots <i_b \leq n$ of size $b$.

     For later use, we next show $\max_{1 \leq i \leq n} \sigma_i^2/b \rightarrow 0$.  Without loss of generality, we suppose $\sigma_1^2 = \max_{1 \leq i \leq n} \sigma_i^2$ and we average the variances of sample means
$(X_{i_1,\mu}+\cdots + X_{i_b,\mu})/\sqrt{b}$ over all ${n-1 \choose b-1}$ subsamples containing $i_1=1$ to obtain
\begin{eqnarray*}
a_n &\equiv &\frac{1}{{n-1 \choose b-1}} \sum_{1 =i_1<i_2<\cdots <i_b \leq n} \left( \frac{1}{b} \sum_{j=1}^b \sigma_{i_j}^2 \right)\\
 &= &\frac{\sigma_1^2}{b} + \frac{b-1}{b(n-1)} \sum_{i=2}^n \sigma^2_i =\frac{\sigma_1^2}{b} \frac{n-b}{n-1} + \frac{n(b-1)}{b(n-1) } v_n, \end{eqnarray*}
where $v_n \equiv n^{-1}\sum_{i=1}^n \sigma_i^2$.  From assumptions, $a_n \to \sigma^2$ and $v_n\to \sigma^2$, so that $ \max_{1 \leq i \leq n} \sigma_i^2/b \equiv  \sigma_1^2/b \leq
 (n-1)^{-1}(n-b)|a_n - n(b-1) v_n/[b(n-1)]|\to 0 $ because $b^{-1}+b/n\to 0$.

To show $Y_b \stackrel{d}{\to} N(0,\sigma^2)$, it suffices to show the characteristic function
\[\E e^{ \imath t Y_b} =  \frac{1}{{n \choose b }} \sum_{1 \leq i_1<i_2<\cdots <i_b \leq n} \E e^{\imath t (X_{i_1,\mu}+\cdots+X_{i_b,\mu})/\sqrt{b}} \rightarrow e^{- t^2 \sigma^2/2}\] converges for given $t\in\mathbb{R}$ (with $\imath=\sqrt{-1}$).
For a given subsample $1 \leq i_1<i_2\cdots <i_b \leq n$, note $\E e^{\imath t (X_{i_1,\mu}+\cdots+X_{i_b,\mu})/\sqrt{b}}= \prod_{j=1}^b e^{\imath t X_{i_j,\mu}/\sqrt{b}}$ by independence and, using  $t^2 \max_{1 \leq i \leq n} \sigma_i^2/b<1$ for all large $n$ along with $|\prod_{i=1}^b w_i -\prod_{i=1}^b z_i |\leq \sum_{i=1}^b |w_i-z_i|$ for complex numbers $w_i,z_i$ with $|w_i|,|z_i|\leq 1$, we may bound
\begin{eqnarray*}
&&\left|\prod_{j=1}^b \E e^{\imath t  X_{i_j,\mu}/\sqrt{b}} - \prod_{j=1}^b  [ 1 -   t^2 \sigma_{i_j}^2/(2b)]
+\prod_{j=1}^b [ 1 -   t^2 \sigma_{i_j}^2/(2b)] -   \prod_{j=1}^b  e^{- t^2 \sigma_{i_j}^2/(2b)} \right|\\&\leq&\sum_{j=1}^b \left|\E  e^{\imath t X_{i_j,\mu}/\sqrt{b}} - [1 - t^2 \sigma_{i_j}^2/(2b)]      \right| + \sum_{j=1}^b
\left| 1 -   t^2 \sigma_{i_j}^2/(2b)  - e^{- t^2 \sigma_{i_j}^2/(2b)}  \right|
\end{eqnarray*}
where further \[\left|\E  e^{\imath t X_{i_j,\mu}/\sqrt{b}} - [1 - t^2 \sigma_{i_j}^2/(2b)]      \right| \leq \E \min\{X_{i_j,\mu}^2/b, |X_{i_j,\mu}|^3/b^{3/2}\}\]
by $|e^{\imath x} - [1 + \imath x + (\imath x)^2/2]|\leq \min\{|x|^2,|x|^3/3!\}$, $x\in\mathbb{R}$, $\E X_{i_j,\mu}=0$ and $\E X_{i_j,\mu}^2=\sigma^2_{i_j}$.
Also,  $| 1 -   t^2 \sigma_{i_j}^2/(2b)  - e^{- t^2 \sigma_{i_j}^2/(2b)}  |\leq [t^2 \sigma_{i_j}^2/(2b)]^2 e^{t^2 \sigma_{i_j}^2/(2b)}$
by $|e^x - 1 -x| \leq x^2 e^{|x|}$, $x\in\mathbb{R}$.   Hence, we may bound
\[
  | \E e^{ \imath t Y_b} -e^{- t^2 \sigma^2/2} | \leq \Delta_{1n}+\Delta_{2n}+\Delta_{3n}
\]
with
\begin{eqnarray*}
\Delta_{1n}& \equiv & \frac{1}{{n \choose b }} \sum_{1 \leq i_1<i_2<\cdots <i_b \leq n} \sum_{j=1}^b \E \min\{X_{i_j,\mu}^2/b, |X_{i_j,\mu}|^3/b^{3/2}\}
\\&=& \frac{1}{n} \sum_{i=1}^n \E \min\{X_{i,\mu}^2, |X_{i,\mu}|^3/\sqrt{b}\},
\end{eqnarray*}
and where
\begin{eqnarray*}
\Delta_{2n} &\equiv&  \frac{1}{{n \choose b }} \sum_{1 \leq i_1<i_2<\cdots <i_b \leq n} \sum_{j=1}^b [t^2 \sigma_{i_j}^2/(2b)]^2 e^{t^2 \sigma_{i_j}^2/(2b)} \\&\leq& t^4 \left(\frac{1}{b}\max_{1 \leq i \leq n} \sigma^2_i\right) \left(e^{ t^2  \max_{1 \leq i \leq n} \sigma^2_i/b}\right)  v_n   \rightarrow 0
\end{eqnarray*}
by $v_n \equiv n^{-1}\sum_{i=1}^n \sigma_i^2 \to \sigma^2$ and $\max_{1 \leq i \leq n} \sigma^2_i/b\to 0$ and where
\[
\Delta_{3n} \equiv  \max_{1 \leq i_1<i_2<\cdots <i_b \leq n}\left| e^{-t^2 \sum_{j=1}^b \sigma^2_{i_j} /(2 b) } - e^{-t^2 \sigma^2/2} \right|\to 0
\]
by continuity and the assumptions.   For a given $\epsilon>0$, we may write
\[
 \Delta_{1n} \leq   \frac{1}{n} \sum_{i=1}^n \E  X_{i,\mu}^2 I(|X_{i,\mu}| > \epsilon \sqrt{b}) + \frac{\epsilon}{n} \sum_{i=1}^n \E  X_{i,\mu}^2 I(|X_{i,\mu}| \leq \epsilon \sqrt{b})
\]
so that $\overline{\lim}_{n\to \infty}  \Delta_{1n} \leq \epsilon \sigma^2$ by  $ n^{-1}\sum_{i=1}^n \E  X_{i,\mu}^2 I(|X_{i,\mu}| > \epsilon \sqrt{b})\to 0$ from assumption and $v_n \equiv n^{-1}\sum_{i=1}^n   \E  X_{i,\mu}^2  \to \sigma^2$.  This shows $\Delta_{1n}\to 0$
and $\E e^{ \imath t Y_b} \rightarrow e^{- t^2 \sigma^2/2}$.

To show convergence of $\Cniid$ and $\Cniidb$ in Mallow's metric $d_2(\cdot,\cdot)$, we define independent bootstrap sample $X_1^*,\ldots,X_b^*$
and a corresponding subsample-like resample $Y_1^*,\ldots,Y_b^*$ from the data $\{X_i\}_{i=1}^n$ as follows.  Let $I_1,\ldots,I_b$ be iid variables, drawn uniformly from $\{1,\ldots,n\}$,
and define $X_j^* \equiv X_{I_j}$ for $j=1,\ldots,b$ and $\bar{X}_b^* \equiv \sum_{j=1}^b X_j^*/b$.  We then set $\bar{Y}_b^* \equiv \sum_{j=1}^b Y_j^*/b$ and $Y_j^* = X_{J_j}$, $j=1,\ldots,b$, using index random variables $J_1,\ldots,J_b$ defined as follows: set $J_1=I_1$ and,
for $i>1$, set $J_i=I_i$ if $I_i \notin \{J_1,\ldots,J_{i-1}\}$ and otherwise choose $J_i$ uniformly from $\{1,\ldots,n\}\setminus \{J_1,\ldots,J_{i-1}\}$.  In this construction, the indices $J_1,\ldots,J_b$ form a (no replacement) subset of $\{1,\ldots,n\}$ and, due to independence between $I_i$ and $\{J_1,\ldots,J_{i-1}\}$ and uniformity, it holds that the resample probability $P_*( I_i=J_i ) = (n+1-i)/n$, $i=1,\ldots,b$, and
the   resampling distribution of $\sqrt{b} (\bar{Y}_b^*-\bar{X}_n)$  matches the subsampling distribution $\Lni$ (i.e., $P_*( \sqrt{b} (\bar{Y}_b^*-\bar{X}_n) \leq x) = \Lni(x)$, $x\in\mathbb{R}$).  Both $\bar{Y}_b^*$ and $\bar{X}_b$ have the same resample mean, $\E_* \bar{Y}_b^* = \E_*\bar{X}_b^*=\bar{X}_n$, and
 both distributions $\Cniid$ and $\Cniidb$ correspond to scaled sums of $k_n$ iid terms distributed, respectively, as either $\sqrt{b} (\bar{Y}_b^*-\bar{X}_n)$ or $\sqrt{b} (\bar{X}_b^*-\bar{X}_n)$. Consequently (cf.~ Lemma~8.7, \citep{BiFr}), we may bound the squared Mallow's distance as
 \[
 [d_2(\Cniid,\Cniidb)]^2 \leq  b \E_*(\bar{Y}^*_b-\bar{X}_b^*)^2 =\frac{1}{b}\sum_{i,j=1}^b \E_*(X_i^* - Y_i^*)(X_j^* - Y_j^*).
 \]
  From the resampling construction, note that  $\E_*[(X_i^* - Y_i^*)(X_j^* - Y_j^*)| I_i=J_i \,\mbox{or}\, I_j=J_i]=0$ for any $1\leq i,j \leq b$.  Also, using conditional probability and uniformity, we have, for any $1 \leq i \neq j \leq b$, that
 \begin{eqnarray*}
 && |\E_*[(X_i^* - Y_i^*)(X_j^* - Y_j^*)| I_i\neq J_i, I_j \neq J_j]|\\
 &\leq&
  \frac{1}{n^2(n-1)(n-2)} \left|\sum_{1\leq i,j,k,m \leq n \atop i\neq k, k \neq m, j \neq m}
  (X_i-X_k)(X_j-X_m)\right|
  \leq  \frac{8}{n} W_n
  \end{eqnarray*}
for $W_n \equiv n(\bar{X}_n-\mu)^2 + n^{-1}\sum_{i=1}^n X_{i,\mu}^2$;  likewise, for $i=1,\ldots,b$,
\[
  |\E_*[(X_i^* - Y_i^*)^2| I_i\neq J_i]| \leq \frac{1}{n(n-1)}  \sum_{1 \leq i \neq j \leq n} (X_i-X_j)^2 \leq \frac{8}{n} \sum_{i=1}^n X_{i,\mu}^2\leq 8 W_n.
\]
Hence, using $ P_*(I_i\neq J_i, I_j \neq J_j) \leq (\min\{i,j\}-1)/n$ from above, we have
\begin{eqnarray*}
 [d_2(\Cniid,\Cniidb)]^2 &\leq  &  \frac{1}{b}\sum_{i=2}^{b }  8 W_n \frac{i-1}{n} \left(1+    2 \sum_{j= i}^{b} \frac{1}{n}\right)\\&\leq &      16 \frac{b}{n}\left(1+\frac{b}{n}\right) W_n =o_p(1)
 \end{eqnarray*}
using that $W_n =O_p(1)$ by $\E |W_n| = 2 v_n \equiv 2n^{-1}\sum_{i=1}^n \sigma_i^2 \to \sigma^2$ and $b/n\to 0$. $\Box$

\end{appendix}

\end{document}